\documentclass[reqno]{amsproc}

\usepackage{amstext,amsmath,amssymb,amsfonts}
\usepackage[latin1]{inputenc}
\usepackage{epsfig}
\usepackage{hyperref}
\usepackage{color}

\usepackage{latexsym}

\DeclareMathOperator{\rank}{rank}

\newtheorem{lemma}{Lemma}
\newtheorem{corollary}{Corollary}

\newtheorem{definition}{Definition}
\newtheorem{theorem}{Theorem}

\newtheorem{proposition}{Proposition}
\newtheorem{example}{Example}

\newcommand{\vspi}{\vspace{0.4cm}}

\newcommand{\bee}{\mathbf{b}}

\newcommand{\bea}{\begin{eqnarray}}
\newcommand{\eea}{\end{eqnarray}}
\newcommand{\beq}{\begin{equation}}
\newcommand{\eeq}{\end{equation}}
\newcommand{\enn}{\nonumber \end{equation}}

 \newcommand{\rk}{{\rm r}\,}

 \newcommand{\cG}{\mathcal{G}}
 \newcommand{\cH}{\mathcal{H}}
 \newcommand{\cV}{\mathcal{V}}
 \newcommand{\cE}{\mathcal{E}}
 \newcommand{\cF}{\mathcal{F}}

 \newcommand{\cB}{\mathcal{B}}
 
 \newcommand{\cR}{\mathcal{R}}

\newcommand{\mT}{\mathfrak{T}}

\newcommand{\bG}{\partial{\mathcal{G}}}
\newcommand{\bV}{{\mathcal{V}}_{\partial}}
\newcommand{\bE}{{\mathcal{E}}_{\partial}}
\newcommand{\bC}{{\mathcal{C}}_{\partial}}

\newcommand{\mG}{\mathfrak{G}}

 \newcommand{\mf}{\mathfrak{f}}

 \newcommand{\sset}{\Subset}

 \newcommand{\inter}{{\rm int}}
 \newcommand{\ext}{{\rm ext}}

\title[Separate half-edged ribbon graphs and tensor graphs]{The expansion of polynomial invariants for $2$-decompositions of generalized graphs}

\author{Remi Cocou Avohou}
\address[R.C.A.]{
International Chair in Mathematical Physics and Applications,
ICMPA-UNESCO Chair, University of Abomey-Calavi, 072BP50, Cotonou, Rep. of Benin, and African Institute for Mathematical Sciences (AIMS-Senegal)}
\email{avohouremicocou@yahoo.fr}

\begin{document}

\maketitle 

\begin{abstract}
The $2$-decomposition for ribbon graphs was introduced in [Annals of Combinatorics 15  (2011), pp 675-706]. We extend this result to  half-edged ribbon graphs and to rank $D$-weakly colored graphs [SIGMA 12 (2016), 030], generalizing therefore the $2$-sums and tensor products of these graphs. Using this extension for the $2$-decompositions, we provide new expansion formulas for  the Bollob\'as Riordan polynomial  for half-edged ribbon graphs and also for the polynomial invariant for weakly colored stranded graphs.
\\

\end{abstract}

\

\begin{center}
Dedicated to Mahouton Norbert Hounkonnou's 60th birthday anniversary
\end{center}

\

\tableofcontents

\section{Introduction}

A graph-theoretic invariant called separability was recently introduced by Cicalese and  Milan${\rm \check{c}}$  \cite{cic}. In general, a graph is called $k-$separable if any two non-adjacent vertices can be separated by the removal of at most $k$ vertices. The $k-$separability turns out to be an important property of a graph used to investigate the computational complexity of several optimization problems for graphs of bounded separability \cite{cic}.
 One of the main
results in \cite{cic} is a decomposition theorem for the 2-separable graphs. 
This decomposition has been extended to ribbon graphs 
by taking into account the cyclic order of the vertices  \cite{stephen-fathom}. 

Defined as a neighborhood of a graph embedded in a surface, a ribbon graph \cite{bollo,joamo} can be decomposed into its $2-$connected components according to the following description: assume that a ribbon graph $\widehat{\cG}$ is $2-$separable and then regard it as arising from a sequence of $2-$sums of a collection of ribbon graphs $\{A_e\}_{e\in\cE}$  with a ribbon graph $\cG=(\cV,\cE)$.  Along the edge $e\in 
\cE$, we ``glue'' the distinguished ribbon graph $A_e$. Strictly speaking, a $2-$sum of two graphs $\cG$ and $\cF$ with distinguished edges $e$ and $f$ respectively in $\cG$ and $\cF$ is defined by identifying $e$ with $f$ and deleting the identified edge \cite{stephen-fathom}. The structure $(\cG,\{H_e\}_{e\in\cE})$, is called a $2-$decomposition for $\widehat{\cG}$ with $H_e=A_e-e$, for all $e \in \cE$.

In \cite{stephen-fathom}, Huggett and Moffat  find a connection between the Bollob\'as Riordan (BR) polynomial of $\widehat{\cG}$ and those of $\cG$ and $H_e$. This result  generalizes Brylawski's  results in \cite{Bryl} which used the universal properties of the Tutte polynomial \cite{tutte}. The goal of the present work is to make one step further and to extend these series of results to new classes of generalized graphs appearing in quantum 
field theory and in theoretical physics \cite{Gurau:2009tz,Ben}. 

Graphs can be generalized as half-edged graphs  (HEGs)  \cite{rca}. 
A half-edge is defined as any edge incident to a unique vertex without forming a loop and a
HEG is a graph together with an incidence relation which associates each half-edge with a unique vertex.  HEGs are the natural class of graphs of quantum field theory \cite{krf},
with half-edges representing field modes with much lower energy than internal processes represented by well formed edges. 
Combining the definitions of HEGs and ribbon graphs, half-edged ribbon graphs (HERGs) arose as a class of graphs encompassing those. HERGs have both an underlying half-edged graph and  ribbon graph structures. Ribbon graphs are also
surfaces with boundary and each boundary component is called
face of the ribbon. In the case of HERGs,  the presence of the half-edges induces two kinds of faces: internal and external faces. The internal faces are components  homeomorphic to $S_1$ and the external are the remaining ones which are homeomorphic to any open segment. Following the external faces, we obtain connected components called 
connected components of the boundary graph. As a new feature, the notion of $2-$decomposition
 for HERGs that we introduce in the present work  distinguishes the treatment 
of the internal faces and the connected components of the boundary graph during the $2-$sum operation.

The Tutte and  Bollob\'as Riordan polynomials have found an extension to HERGs by including extra variables: one for keeping track  of the number of connected components of the boundary graph and another for the number of half-edges. Using the bijection between the states of $\widehat{\cG}$ and those of $\cG$ and $H_e$, we find
the expression for the number of internal faces, number of connected
components of the boundary graph and half-edges of a state of $\widehat{\cG}$ 
in terms of those of the states of $\cG$ and $H_e$. This task remains complex
because, first, the definition by Hugget and Moffat must be modified to take into
account the presence of the half-edges. Then, we must deal with the fact that the internal faces or connected components of the boundary graph of a state in $\widehat{\cG}$ may be generated by different types of components in the states in $\cG$ and the graphs $H_e$. To tackle this issue, we identify 
a matrix $\epsilon$ which captures the subtlety of  the $2-$decomposition of HERGs. 
This matrix called $\epsilon(f,x_{e})$ is labelled by rows indexed by the internal faces of a state in $\cG$ and columns indexed by the points $x_{e}$ that define the meeting places of the graphs $H_e$ on $\cG$ in the construction of $\widehat{\cG}$. The important quantity to master the decomposition of HERGs is the rank of $\epsilon$.
As a consequence, we relate
the multivariate polynomial invariant of $\widehat{\cG}$ and of those
of  $\cG$ and of $H_e$.

 The level of difficulty increases when we seek for such relations for the class of graphs called rank $D$ weakly-colored (w-colored) stranded graphs \cite{Avo16}.
Such graphs are called  stranded graphs because they are made with stranded vertices which are chord diagrams and stranded edges which are collections of segments. Gurau introduced in \cite{Gurau:2009tzb} a coloring on them and proved that  they are dual to simplicial pseudo-manifolds in any dimension $D$. This duality was the stimuli for particular quantum gravity models claiming that the geometry of spacetime
at high energy is simplicial.  In that sense, stranded graphs represent quantum (discrete)
spaces. 
In a colored graph, the vertices are called $0-$cells, lines or edges $1-$cells and the faces $2-$cells. A $p-$cell or $p-$bubble is defined as a connected subgraph made only of lines of $p$ chosen colors. Once we impose to the vertices to have a fixed coordination $D$, we obtain a specific stranded graph called colored tensor graph \footnote{Notice that this notion of tensor graph is technically different from the notion of tensor product of graphs defined as a particular 2-decomposition.}. The coordination of the vertices in a colored tensor graph gets modified if we perform a contraction of an edge. Allowing such a contraction enlarges the class of graphs from the colored tensor graphs to what is called  weakly-colored (w-colored)  graphs. The BR polynomial invariant has a natural extension from ribbon graphs to rank $D$ w-colored graphs \cite{rca, Avo16}.
In the second part of this work, we introduce the 2-decomposition for rank $D$ w-colored 
graphs and establish few properties of it. The expression of the generalized invariant
on the graph $\widehat{\cG}$ in terms of the invariants of the corresponding graphs
$\cG$ and $H_e$ is much more involved.

In this paper, section \ref{sect:preliminaries} reviews the ribbon graphs, half-edged ribbon graphs and the BR polynomials on such graphs.  We show in Section \ref{sect:EXpBR}, how to compute the BR polynomial of a graph $\widehat{\cG}=(\cG,\{H_e\}_{e\in\cE})$, by introducing two matrices which capture the details generated by that $2-$decomposition.
 Two particular cases are studied: the case where each graph $H_e$ is embedded in the neighborhood of $e$ in the embedded graph $\cG$ and the general case. Theorems \ref{theo:prod} and \ref{theo:import}  establish the main results of this section. 
We investigate the 2-decomposition extended to the class of weakly colored graphs in section \ref{sect:weaklycoloreddec} and report already a preliminary result
 of an expansion of the invariant of the weakly colored graphs. 
The complete expansion in terms of the invariants of $\cG$ and of
$H_e$ is solved for a particular class of weakly colored graphs. 

\section{Preliminaries}
\label{sect:preliminaries}
In this paper, we briefly review some essential concepts on ribbon graphs, HERGs and the Bollob\'as Riordan polynomial.

A ribbon graph $\cG=(\cV,\cE)$ is a surface with boundary where the vertices are represented by a set of disks and the edges are represented by ribbons (rectangular disks) \cite{bollo}. A spanning subgraph $s=(\cV,\cE')$ of $\cG$, where $\cE'\subseteq \cE$, is called a state of $\cG$ and we denote by $S(\cG)$ the set of states of $\cG$. Let $v(s)$, $e(s)$, $k(s)$, $r(s)$, $n(s)$ and $\partial(s)$ be respectively the number of vertices, edges, connected components, rank, nullity and boundary components of $s$. Besides the parameters of this graph, there is $t(s)$, which records the orientability of an embedded graph $s$. By definition $t(s) = 0$ if $s$ is orientable and $t(s) = 1$ otherwise.

There are some graphs operations which we now describe

\begin{definition}[Deletion and contraction \cite{bollo}] 
\label{def:cdelrib}

Let $\cG$ be a ribbon graph and $e$ one of its edges. 

$\bullet$ We call $\cG-e$ the ribbon graph obtained from 
$\cG$ by deleting $e$ and keeping the end vertices 
as closed discs. 

$\bullet$ If $e$ is not a loop and is positive, consider its
end vertices $v_1$ and $v_2$. The graph $\cG/e$  obtained by contracting $e$ is defined from $\cG$ by replacing $e$, 
$v_1$ and $v_2$ by a single vertex disc $e\cup v_1 \cup v_2$.
If $e$ is a negative non-loop, then untwist it (by flipping
one of its incident vertex) and  contract.

$\bullet$ If $e$ is a trivial twisted loop, 
contraction is deletion: $\cG-e = \cG/e$. 
The contraction of a trivial untwisted loop  $e$
is the deletion of the loop and the addition of a new 
connected component vertex $v_0$ to 
the graph  $\cG-e$. We write $\cG/e = (\cG-e)\sqcup \{v_0\}$.

\end{definition}

The Bollob\'as-Riordan polynomial $R(\cG; x,y,z,w) \in Z[x,y,z,w]/w^2 - w$ for ribbon graphs is defined as the state sum:
\bea\label{eq:bolo}
R(\cG; x, y, z, w) = \sum_{s\in S(\cG)}
(x - 1)^{r(\cG)-r(s)} y^{n(s)} z^{k(s)- \partial(s)+n(s)}w^{t(s)} .
\eea
In the following we set $w=1$ and use an abuse of notation $e\in s$ rather than $e\in \cE(s)$.
Under these assumptions, the multivariate Bollob\'as-Riordan polynomial \cite{vigne} is
\bea\label{eq:bolomuti}
Z(\cG; a, b, c) = \sum_{s\in S(\cG)} a^{k(s)}\Big(\prod_{e\in s}b_e\Big)c^{\partial(s)},
\eea
where $a$ and $c$ are indeterminates, and $b := \{b_e /e \in \cE\}$ is a set of indeterminates indexed by $\cE$. If we set all the variables $b_e = b$ in \eqref{eq:bolomuti}, then using equation \eqref{eq:bolo}, we obtain 
\bea
R(G; x, y, z) = (x-1)^{-k(\cG)}(yz)^{-v(\cG)}Z(\cG; (x-1)yz^2, yz, z^{-1}). 
\eea

Let us now discuss a polynomial invariant for ribbon graphs with half-edges or half-ribbons. Note that the HERGs was originally studied in \cite{krf} where a half-ribbon edge (or simply half-ribbon, denoted henceforth HR) is a ribbon incident to a unique vertex by a unique segment and without forming loops.  Half-ribbons allowed another graph operation called the cutting of an edge.
Cutting an edge $e$ in a ribbon graph $\cG$, means that we remove
$e$ and we let two HRs attached at the end vertices of $e$. If $e$ is a loop, the two HRs are on the same vertex.

\begin{figure}[h]
 \centering
     \begin{minipage}[t]{.8\textwidth}
      \centering
\includegraphics[angle=0, width=6.5cm, height=2cm]{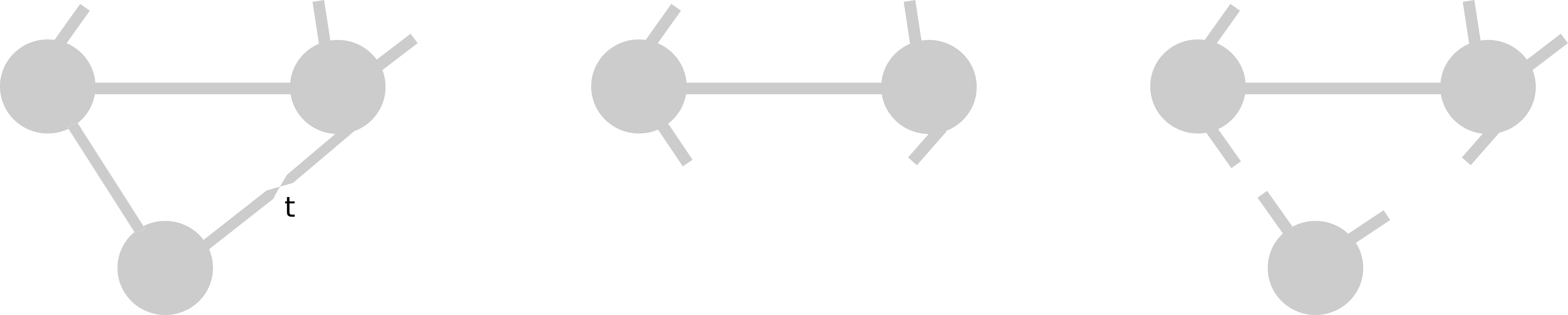}
\vspace{0.3cm}
\caption{ {\small A ribbon graph with HRs $\cG_{ \mf^{0}}$ together with 
a c-subgraph $A_{ \mf^{0}_A}$ and a spanning c-subgraph $s_{\mf^0_s}$ }}
\label{fig:ribgraph}
\end{minipage}
\put(-235,-10){$\cG_{ \mf^{0}}$}
\put(-175,-10){$A_{ \mf^{0}_A}$}
\put(-97,-10){$s_{\mf^0_s}$}
\end{figure}

$\bullet$ A ribbon graph $\cG$ with HRs or a HERG is defined as a ribbon graph $\cG(\cV,\cE)$ with a set $\mf=\mf^1\cup\mf^0$, where $\mf^1$ is the set of HRs obtained from the cut of all edges of $\cG$, and $\mf^0$ is the set of additional HRs together with a relation which associates  with each additional HR a unique vertex. A ribbon graph $\cG(\cV,\cE)$ with the set $\mf^0$ of additional HRs  is denoted $\cG_{\mf^0}(\cV,\cE)$. An illustration is given in Figure \ref{fig:ribgraph}. The subgraphs of $\cG_{\mf^0}(\cV,\cE)$ are obtained by using the ``cutting'' operation to replace the usual operation of deletion.

$\bullet$ A c-subgraph $A_{\mf^0_A}$ of $\cG_{\mf^0}(\cV,\cE)$ is defined as a HERG $A_{\mf^0_A}(\cV_A,\cE_A)$ satisfying $\cV_A\subseteq \cV$ and $\cE_A\subseteq\cE$ such that the incidence relation between edges and vertices is respected. We now denote by $\cE_A'$ the set of edges incident to the vertices of $A$ and not contained in $\cE_A$. The set of HRs of $A_{\mf^0_A}$  is $\mf^0_A = \mf^{0;0}_A \cup \mf^{0;1}_A(\cE_A)$ with $\mf^{0;0}_A \subseteq \mf^0$ and $\mf^{0;1}_A (\cE_A)\subseteq \mf^1$, where $\mf^{0;1}_A (\cE_A)$ is the set of HRs obtained by cutting all edges in $\cE_A'$ and incident to the vertices of $A_{\mf^0_A}$. 
We denote $A_{\mf^0_A} \subseteq \cG_{\mf^0}$. See an example of c-subgraph $A_{\mf^0_A}$ in Figure \ref{fig:ribgraph}.

$\bullet$ A spanning c-subgraph $s_{\mf^0_s}$ or a state of $\cG_{\mf^0}(\cV,\cE)$ is defined as   a c-subgraph $s_{\mf^0_s}(\cV_s,\cE_s)$ of $\cG_{\mf^0}$ with all vertices  and all additional HRs of $\cG_{\mf^0}$. Then $\cE_s\subseteq \cE$ and $\cV_s = \cV$, $\mf^0_s = \mf^{0} \cup \mf^{0;1}_s(\cE_s)$. 
We use the notation $s_{\mf^0_s}\sset \cG_{\mf^0}$ and denote by $S(\cG_{\mf^0})$ the set of states of $\cG_{\mf^0}$. 
(See $s_{ \mf^{0}}$ in Figure \ref{fig:ribgraph}.)

The  states or spanning c-subgraphs have crucial importance in this framework since they are involved in the state sum of  the BR polynomial defined on HERGs. 

There are two kind of boundary on a HERG which deserve to be analyzed:  the boundary faces following the contour of the HRs and the initial ones which follow the boundary of well-formed edges. More precisely let us consider  a ribbon graph with HRs $\cG_{\mf^0}(\cV,\cE)$.

$\bullet$ A closed or internal face is defined as a boundary component of a ribbon graph  which never passes through any free segment of the additional HRs.  We denote by $\cF_{\inter}(\cG_{\mf^0}(\cV,\cE))$ the set of closed faces of $\cG_{\mf^0}(\cV,\cE)$. 

$\bullet$ A boundary  component  obtained by leaving an external point of some HR rejoining another external point is called an open or external face.  We denote by  $\cF_{\ext}(\cG_{\mf^0}(\cV,\cE))$, the set of open faces of $\cG_{\mf^0}(\cV,\cE)$. 

$\bullet$ The set of faces $\cF(\cG_{\mf^0}(\cV,\cE))$ of $\cG_{\mf^0}(\cV,\cE)$ is defined by 
$\cF_{\inter}(\cG_{\mf^0}(\cV,\cE)) \cup \cF_{\ext}(\cG_{\mf^0}(\cV,\cE))$. If  $\cF_{\ext}(\cG_{\mf^0}(\cV,\cE))\neq \emptyset$ i.e. $\mf^0\neq \emptyset$ then $\cG_{\mf^0}(\cV,\cE)$ is said to be open. Otherwise it is closed .

\begin{figure}[h]
 \centering
     \begin{minipage}[t]{.8\textwidth}
      \centering
\includegraphics[angle=0, width=4cm, height=2cm]{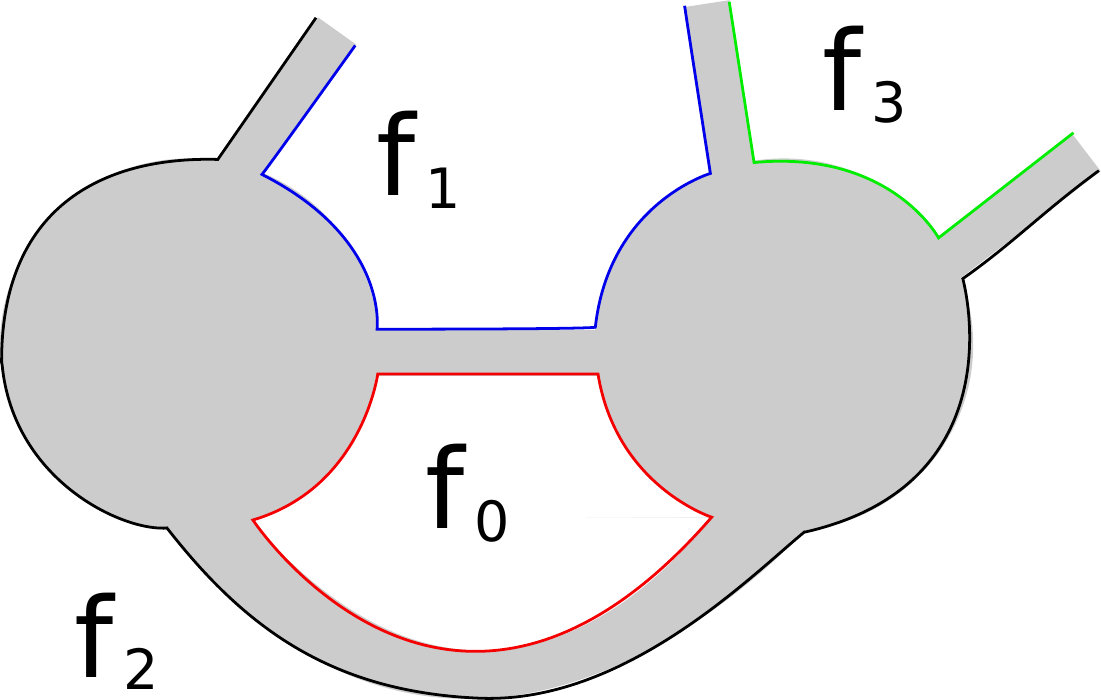}
\vspace{0.3cm}
\caption{ {\small A HERG $\cG_{\mf^0}$ with $\cF_{\inter}(\cG_{\mf^0})=\{f_0\}$, and $\cF_{\ext}(\cG_{\mf^0})=\{f_1,f_2,f_3\}$ }}
\label{fig:faces}
\end{minipage}
\end{figure}

An illustration for closed and open faces is given in Figure \ref{fig:faces}.

The boundary $\partial\cG_{\mf^0}$ of a ribbon graph $\cG_{\mf^0}(\cV,\cE)$ is a simple graph $\partial\cG_{\mf^0}(\bV,\bE)$ such that $\bV$ is one-to-one with $\mf^0$ and $\bE$ is one-to-one with $\cF_{\ext}(\cG_{\mf^0}(\cV,\cE))$. Then the boundary graph of a closed ribbon graph is empty. By construction, the boundary graph $\partial\cG_{\mf^0}$, is obtained by inserting a vertex of valence or degree two at each HR, the edges of $\partial\cG_{\mf^0}$ are nothing but the external faces of $\cG_{\mf^0}$.
\begin{figure}[h]
 \centering
     \begin{minipage}[t]{.8\textwidth}
      \centering
\includegraphics[angle=0, width=2cm, height=1cm]{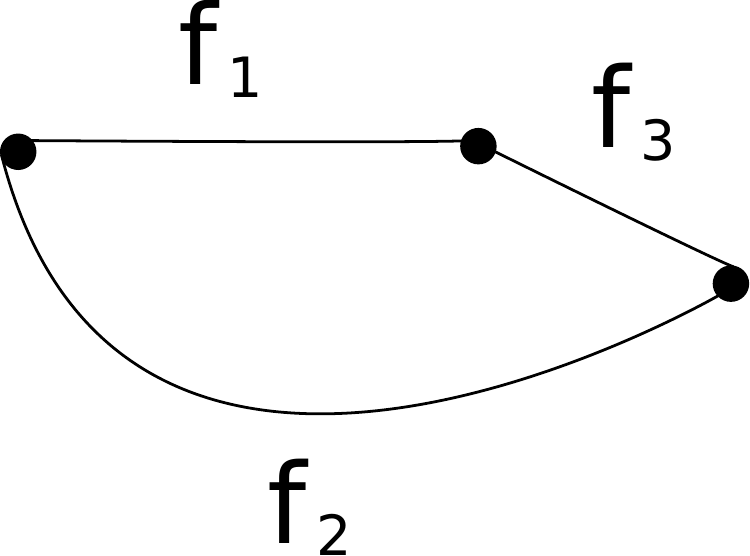}
\vspace{0.3cm}
\caption{ {\small The boundary graph associated to
the HERG in Figure \ref{fig:faces} }}
\label{fig:boundary}
\end{minipage}
\end{figure}
The operations of edge contraction and deletion for HERGs keep their meaning as in Definition \ref{def:cdelrib}. 
\begin{definition}[BR polynomial for HERGs]
Let $\cG_{\mf^0}(\cV,\cE)$ be a HERG. We define the ribbon graph polynomial of $\cG_{\mf^0}$ to be
\bea\label{eq:boloflaede}
\cR_{\cG_{\mf^0}}(x,y,z,w,t)
=\sum_{s\in S(\cG_{\mf^0})} (x-1)^{\rk(\cG_{\mf^0})-\rk(s)}(y-1)^{n(s)}
z^{k(s)-F_{\inter}(s)+n(s)} \, w^{C_\partial(s)}\, t^{f(s)},
\label{brfla}
\eea
where $C_\partial(s)= |\bC(s)|$ is the number of connected components of the boundary of $s$, $f(s)$, the number of half-edges and $F_{\inter}(s)= |\cF_{\inter}(s)|$.
\end{definition}

This definition gives a polynomial $\cR$ \eqref{eq:boloflaede} which is a generalization of the BR polynomial $R$ \eqref{eq:bolo} from ribbon graphs to HERGs. The multivariate Bollob\'as-Riordan polynomial for HERGs is
\bea\label{eq:boloflamuti}
Z_{\cG_{\mf^0}}(a, b, c, d, l) = \sum_{s\in S(\cG_{\mf^0})} a^{k(s)}\Big(\prod_{e\in s}b_e\Big)c^{F_{\inter}(s)}d^{\partial(s)}l^{f(s)},
\eea
where $a$, $c$ and $d$ are indeterminates, and $b := \{b_e /e \in \cE\}$ is a set of indeterminates indexed by $\cE$. If we set all the variables $b_e = b$ in \eqref{eq:boloflamuti}, then using equation \eqref{eq:boloflaede}, we obtain 
\bea\label{eq:boloflaandmult}
\cR_{\cG_{\mf^0}}(x, y, z,w,t) = (x-1)^{-k(\cG_{\mf^0})}(yz)^{-v(\cG_{\mf^0})}Z_{\cG_{\mf^0}}((x-1)yz^2, yz, z^{-1},w,t).
\eea
\section{Expansion for the Bollob\'as-Riordan polynomial on  half-edged ribbon graphs}
\label{sect:EXpBR}
\begin{figure}[h]
 \centering
     \begin{minipage}[t]{.8\textwidth}
      \centering
\includegraphics[angle=0, width=7cm, height=2.5cm]{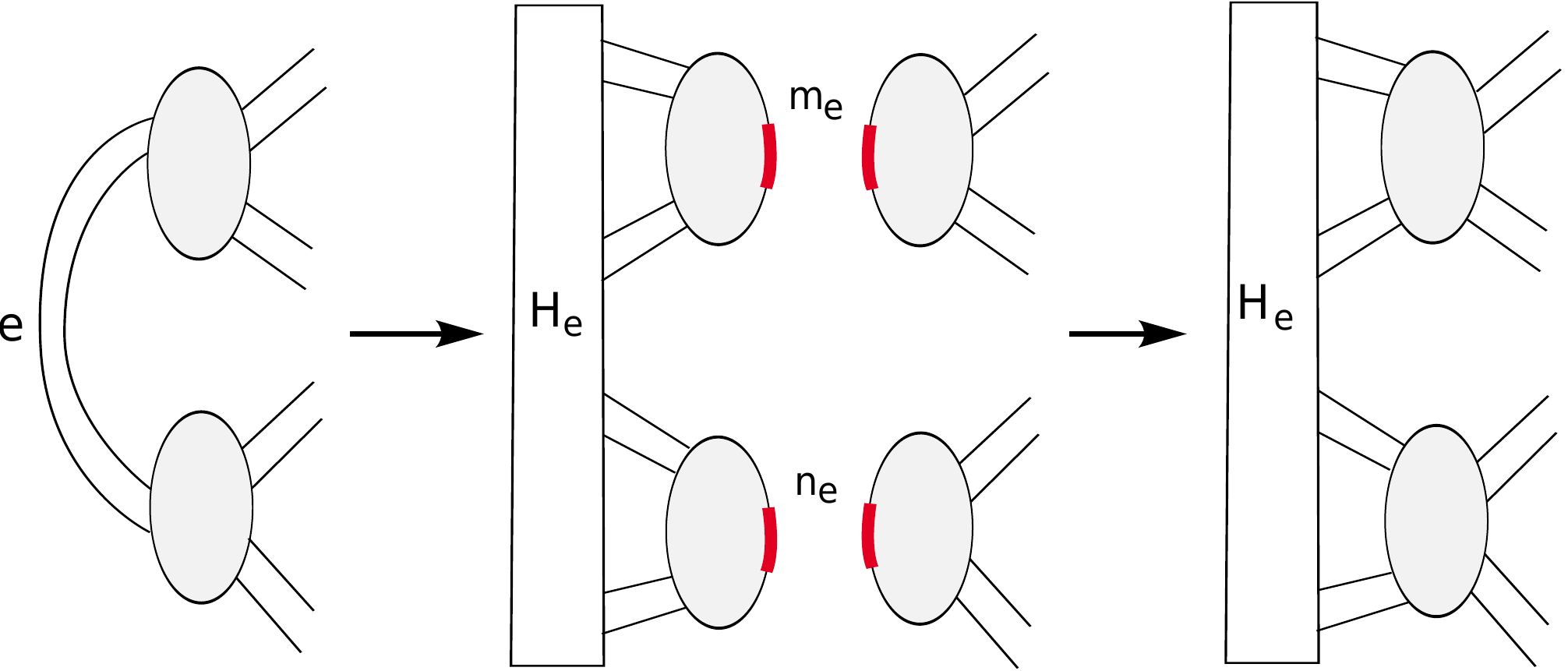}
\caption{ {\small  Replace the edge $e$ by $H_e$ to obtain $\widehat{\cG}_{\mf^0}$}}
\label{fig:gluing}
\end{minipage}
\end{figure}

The formation of the HERG $\widehat{\cG}_{\mf^0}$ from its $2-$decomposition $(\cG_{\mf^0}, \{H_e\}_{e\in\cE})$ is obtained by replacing each ribbon edge $e$ by $H_e$. In fact from the $2$-decomposition $(\cG_{\mf^0},\{H_e\}_{e\in\cE})$ locally at $e=(u_e,w_e)$ of the template $\cG_{\mf^0}$, the graph $\widehat{\cG}_{\mf^0}$ is constructed by identifying the arcs $m_e$ and $n_e$ of the vertices $u_e$ and $w_e$ of $H_e$ with the corresponding arcs $m_e$ and $n_e$ on the vertices $u_e$ and $w_e$ on the template $\cG_{\mf^0}-e$. We use the same notations $n_e$ and $m_e$ according to the identification of the arcs. An illustration is given in Figure \ref{fig:gluing}, where the arcs $m_e$ and $n_e$ are shown in red. 
\subsection{The $2-$decomposition of HERGs} We give in this subsection a natural extension of the $2-$decomposition of graphs, from ribbon graphs to HERGs. Given a HERG $\cG_{\mf^0}=(\cV,\cE)$, we may wonder how to evaluate the number of its  internal faces or connected components of the boundary graph using its $2-$decomposition. 
\begin{definition}
\label{def:graphhat}
Let $\cG_{\mf^0}=(\cV,\cE)$ be a HERG and $\{A_e\}_{e\in\cE}$ be a set of half-edged ribbon graphs each of which has a specific non-loop edge distinguished. For each $e\in \cE$ take the $2-$sum $\cG\oplus_2 A_e$, along the edge $e$ and the distinguished edge in $A_e$, to obtain the half-edged ribbon graph $\widehat{\cG}_{\mf^0}$. For each $e\in\cE$ we define $H_e=A_e- \{e\}$. We will call the structure $(\cG_{\mf^0},\{H_e\}_{e\in\cE})$ a $2-$decomposition for $\widehat{\cG}_{\mf^0}$.
\end{definition}
We can construct the graphs $A_e$ by performing the reverse operation. That is we look $\widehat{\cG}_{\mf^0}$ as a $2-$decomposition $(\cG_{\mf^0},\{H_e\}_{e\in\cE})$ where we identify two distinguished vertices $u_e$ and $w_e$ of each $H_e$ with the corresponding end vertices of $e$ in $\cG_{\mf^0}$. From each $H_e$ we can define a graph $A_e=H_e\cup  \{e\}$. 

The set of states  $S(H_e)$ of each $H_e$ is partitioned into two subsets: $S^{1}(H_e)$ consists of all states in $S(H_e)$ in which $u_e$ and $w_e$ lie in the same connected component, and $S^{2}(H_e)$ consists of all states in $S(H_e)$ in which $u_e$ and $w_e$ lie in different connected component. A state $\hat{s}\in\widehat{\cG}_{\mf^0}$ is obtained by replacing the edges $e$ in a state $s\in S(\cG_{\mf^0})$ with elements of $S^{1}(H_e)$, and the edges $f$ which are not in $s$ by elements of $S^{2}(H_f)$. 
\begin{lemma}\label{lemma:connectedhats}
If  a state $\hat{s}$ of the embedded graph $\widehat{\cG}_{\mf^0}$ is decomposed into states $s\in S(\cG_{\mf^0})$, $s_e\in S^1(H_e)\cup S^2(H_e)$, $e\in\cE$ in the decomposition above, then
\bea\label{rel14}
k(\hat{s})&=&\sum_{e\in\cE} k(s_e) - |\{s_e\in S^1(H_e)\}| - 2|\{s_e\in S^2(H_e)\}| + k(s),
\cr
F_{\inter}(\hat{s}) +  C_\partial(\hat{s}) &=& \sum_{e\in\cE}F_{\inter}(s_e) - |\{s_e\in S^1(H_e)\}| - 2|\{s_e\in S^2(H_e)\}| +  F_{\inter}(s) \cr&+& \sum_{e\in\cE}C_\partial(s_e)+C_\partial(s).
\eea
\end{lemma} 
\proof
The proof of this lemma will follow the one of Lemma 3 in \cite{stephen-fathom} by considering the underling ribbon graphs  $\tilde{\hat{s}}$ and $\tilde{s}_e$ associated respectively to the half edged ribbon graphs $\hat{s}$ and $s_e$ and using the equalities
\bea
\partial(\tilde{\hat{s}}) = F_{\inter}(\hat{s}) +  C_\partial(\hat{s}),\quad \partial(\tilde{s}_e) = F_{\inter}(s_e) + C_\partial(s_e).
\eea
\qed

We are now interested in a separate formula relating $F_{\inter}(\hat{s}) $, $F_{\inter}(s_e) $ and $F_{\inter}(s)$ and a formula relating $C_\partial(\hat{s})$, $C_\partial(s_e)$ and $C_\partial(s)$.
This leads us to consider different cases. The case where a connected component of the boundary of $s$ corresponds to a connected component of the boundary of $\hat{s}$ and the case where a closed face of $s$ corresponds to a connected component of the boundary of $\hat{s}$. We give another definition of the $2-$decomposition in the following subsection in order to overcome this issue.

\subsection{The BR polynomial for HERGs embedding in a neighbourhood}
\label{sect:bropen}
 
In this subsection we assume that the graph $\widehat{\cG}_{\mf^0} = (\cG_{\mf^0},\{H_e\}_{e\in\cE})$ is embedded graph where $\cG$ is embedded and each graph $H_e$ is embedded in the neighborhood of $e$. 

\begin{figure}[h]
 \centering
     \begin{minipage}[t]{.8\textwidth}
      \centering
\includegraphics[angle=0, width=12cm, height=3cm]{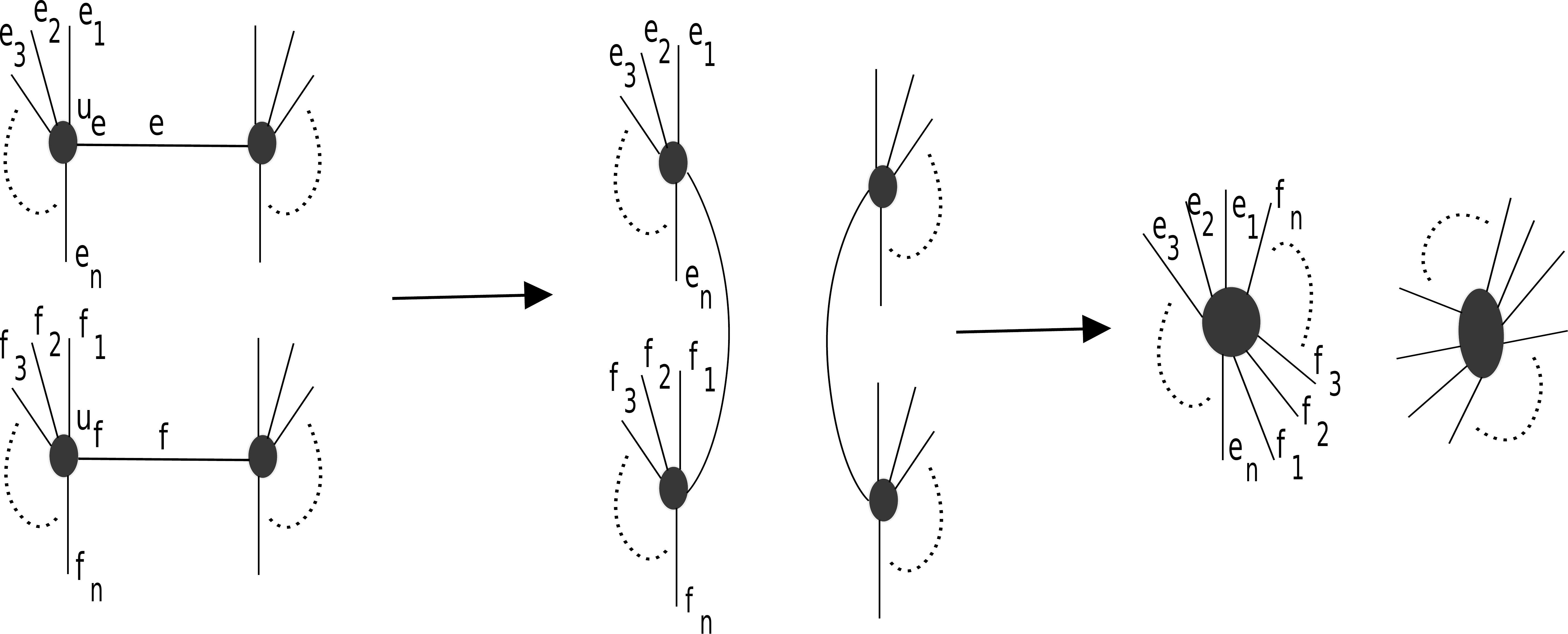}
\caption{ {\small The process of identification of $e$ and $f$ (on the left),  through the introduction  of an edge between  $u_e$ and $u_f$ (midle) and the $2-$sum obtained after contraction (on the right) }}
\label{fig:twosepcolor}
\end{minipage}
\end{figure}

Let us remember the $2-$decomposition in order to evaluate the BR polynomial on  HERGs. Consider two HERGs $\cG_{\mf^0}$ and $\cF_{\mf'^0}$ with distinguished edges $e\in E(\cG_{\mf^0})$ and $f\in E(\cF_{\mf'^0})$. The $2-$sum $\cG_{\mf^0}\oplus_2\cF_{\mf'^0}$ is defined by identifying $e$  with $f$ and deleting the identified edge. We introduced here another way to perform this sum which will be generalized on tensor graphs. In the process of identification of $e$ and $f$, assume that a vertex $u_e\in V(\cG_{\mf^0})$ is identified with $u_f\in V(\cF_{\mf'^0})$. We introduce an edge between $u_e$ and $u_f$ such that the end points of this edge coincide with the end points of the identified edges. We contract this edge and obtain that the cyclic order around the new vertex will be  $\{e_1,\cdots,e_n,f_1,\cdots,f_m\}$ if the cyclic order around $u_e$ and $u_f$ are $\{e,e_1,\cdots,e_n\}$ and $\{f,f_1,\cdots,f_m\}$ according to a choice of orientations on these ribbon graphs. An illustration is given in Figure \ref{fig:twosepcolor}.

\begin{figure}[h]
 \centering
     \begin{minipage}[t]{.8\textwidth}
      \centering
\includegraphics[angle=0, width=6.5cm, height=3cm]{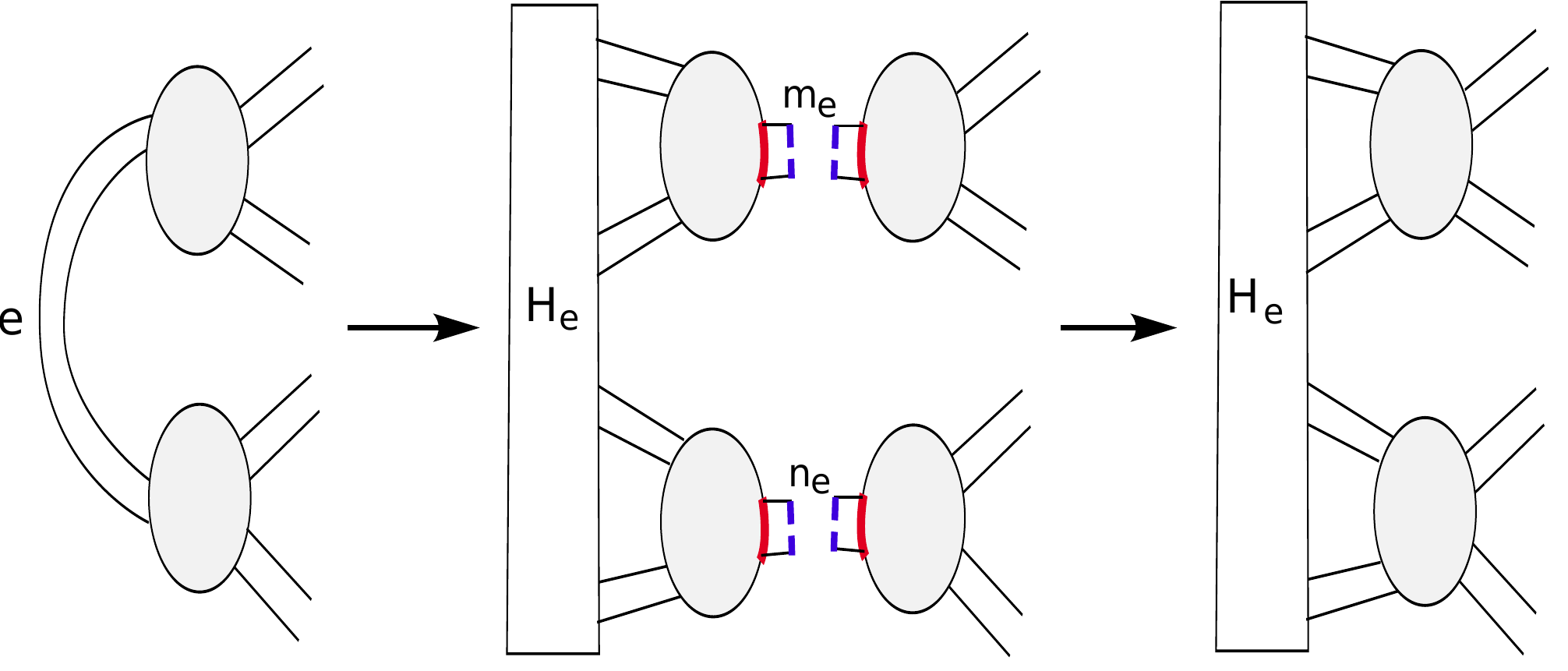}
\caption{ {\small Replace the edge $e$ by $H_e$ to obtain $\widehat{\cG}_{\mf^0}$}}
\label{fig:gluingf}
\end{minipage}
\end{figure}

\begin{definition}
\label{def:graphhatclass2}
Let $\cG_{\mf^0}=(\cV,\cE)$ be a HERG and $\{A_e\}_{e\in\cE}$ be a set of half-edged ribbon graphs each of which has a specific non-loop edge distinguished. For each $e\in \cE$ take the $2-$sum $\cG_{\mf^0}\oplus_2 A_e$, along the edge $e$ and the distinguished edge in $A_e$ as introduced above, to obtain the half-edged ribbon graph $\widehat{\cG}_{\mf^0}$. For each $e\in\cE$ let us define $H_e=A_e\vee \{e\}$. We will call the structure $(\cG_{\mf^0},\{H_e\}_{e\in\cE})$ a $2-$decomposition for $\widehat{\cG}_{\mf^0}$. 
\end{definition} 
The construction of $\widehat{\cG}_{\mf^0}$ from the $2-$decomposition $(\cG_{\mf^0},\{H_e\}_{e\in\cE})$ in Definition \ref{def:graphhatclass2},  is obtained by the identification of the segments called again arcs $m_e$ and $n_e$ of the vertices $u_e$ and $w_e$ in $H_e$ with their correspondence in $\cG_{\mf^0}\vee e$. More precisely these arcs are identified on the two half edges generated by the cutting of $e$ in each graph $A_e$. We give an illustration in Figure \ref{fig:gluingf}, where the arcs $m_e$ and $n_e$ are now the segments in blue. 

Let us make a comparison of the definitions \ref{def:graphhat} and \ref{def:graphhatclass2}.
From Definition  \ref{def:graphhatclass2}, each $H_e$ has at least two half-edges attached to the points $m_e$ and $n_e$. Definition \ref{def:graphhat} shows that $\widehat{\cG}_{\mf^0}$ is obtained by identifying the distinguished vertices $u_e$ and $w_e$,  for $e=(u_e,w_e)$, of each graph $H_e$ with their correspondence in $\cG_{\mf^0}$. Then we may not have some half-edges attached to $m_e$ and $n_e$.

Let $a_e$ and $a'_e$ be the two endpoints of the arcs $m_e$ and $b_e$ and $b'_e$ the endpoints of the arcs $n_e$. Some of these points may belong to the same internal face or connected component of the boundary graph. These points $a_e$, $a'_e$, $b_e$ and $b'_e$ induce points on the boundary of $\cG_{\mf^0}$, $\widehat{\cG}_{\mf^0}$ and $H_e$, and then on the states $s\sset \cG_{\mf^0}$, $\hat{s}\sset \widehat{\cG}_{\mf^0}$ and $s_e\sset H_e$.

Consider one of such points, say, $x$ in an internal face $f$ of a state $s\in S(\cG)$. This point has a correspondent point $x_e$ in each of the graphs $s_e$ such that $f$ will pass by $e$ in $\hat{s}$. We can define a matrix $\epsilon$ which columns are indexed by the points $x_e$ and rows  indexed by  the internal faces of $s$. The elements of this matrix are given by: $\epsilon(f,x_e)=1$ if $x_e$ corresponds to a point $x$ in $f$ and belongs to a connected component of the boundary graph in $s\cup e$ and $\epsilon(f,x_e)=0$ otherwise. To the matrix $\epsilon$ we can associate a sub-matrix $\tilde{\epsilon}$ obtained by removing from $\epsilon$, all the zero and collinear column vectors except one. It is clear that $\tilde{\epsilon}$ is a square matrix. We denote by $\{\tilde x_e\}_e$ the set of the remaining points indexing the columns of $\tilde{\epsilon}$ and $\tau$ the number of internal faces in $s$ which become open in $\hat{s}$. The following result is straightforward.
\begin{lemma}
Let us consider the $2-$decomposition $\hat{s}=(s,\{s\}_e)$ where $\hat{s}\in S(\widehat{\cG})$, $s\in S(\cG)$, $s_e\in S(H_e)$ and $\widehat{\cG}=(\cG,\{H\}_{e})$ and $\epsilon$ as defined above.
We have
\bea
\tau=\rank(\epsilon)=\rank(\tilde{\epsilon}).
\eea
\end{lemma}
We want to express $\tau$ as a sum over e in s. Consider another matrix $\sigma$ whose rows are indexed by the points $\{x_e\}_e$ and  columns indexed   by the edges $e$. We say that $x_g\in e$ if and only if $x_g$ has a correspondent point $x\in e$ and $x\in e=(m_e,n_e)$ means that $x$ is one of the points $a_e$, $a'_e$, $b_e$ and $b'_e$. Each matrix element is given by: $\sigma(x_g,e)=1$ if $x_g\in e$ and $x_g$  belongs to a connected component of the boundary graph in $s\cup e$ and $\sigma(x_e,e)=0$ otherwise. We can deduce the matrix $\tilde\sigma$ indexed by the elements $\{\tilde x_e\}_e$ discussed above. 
\begin{figure}[h]
 \centering
     \begin{minipage}[t]{.8\textwidth}
      \centering
\includegraphics[angle=0, width=10cm, height=5cm]{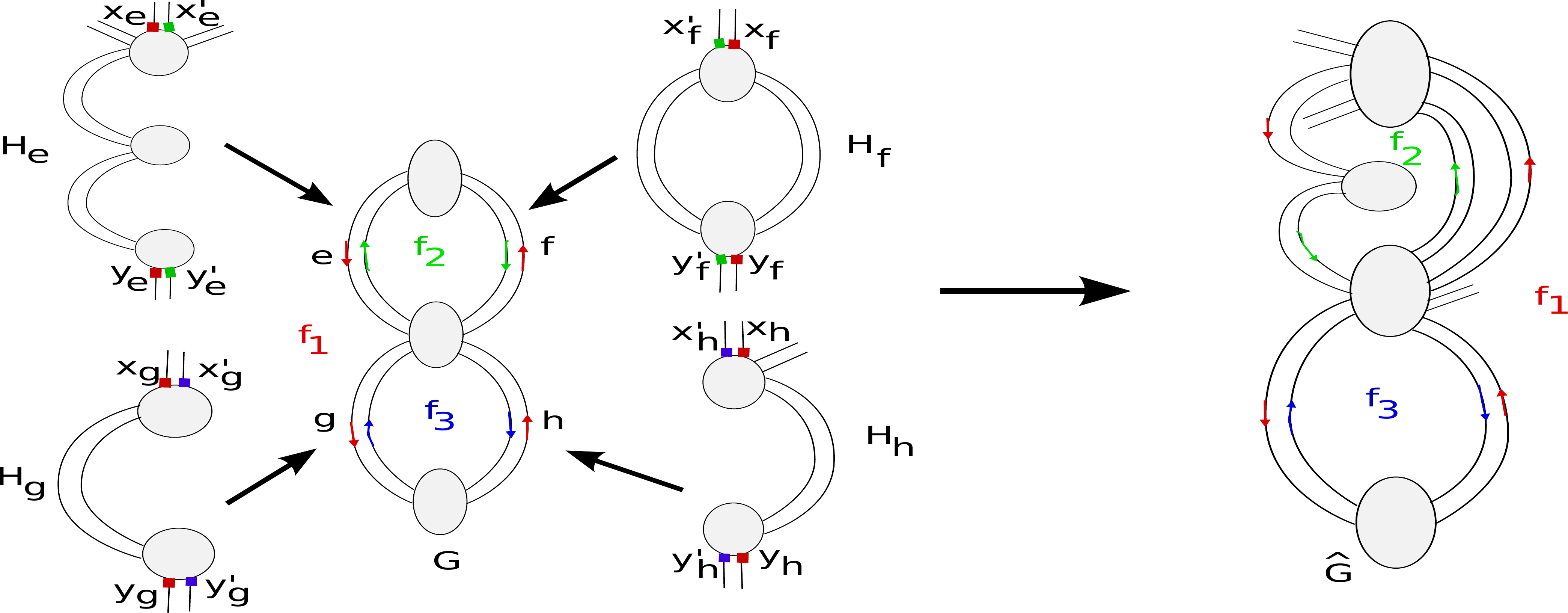}
\caption{ {\small The $2-$decomposition of $\widehat{G}=(G,\{H_e,H_f,H_g,H_h\})$ with the meeting points colored according to the faces of $G$.}}
\label{fig:gluingfpoints}
\end{minipage}
\end{figure}

In Figure \ref{fig:gluingfpoints}, we give the details of the $2-$decomposition of $\hat{G}=(G,\{H_e,H_f,H_g,H_h\})$ in order to compute the matrices $\epsilon$ and $\sigma$. As discussed earlier,  the set of the points is  $$\{x_e,x'_e,x_f,x'_f,x_g,x'_g,x_h,x'_h,y_e,y'_e,y_f,y'_f,y_g,y'_g,y_h,y'_h\},$$ the set of faces, $\{f_1,f_2,f_3\}$ and the edges set $\{e,f,g,h\}$. The columns vectors  in $\epsilon$ indexed by $x_e$, $x'_e$, $x_f$, $x'_f$, $x_g$, $x'_g$, $x_h$ and $x'_h$ are respectively the same with the columns vectors indexed  by $y_e$, $y'_e$, $y_f$, $y'_f$, $y_g$, $y'_g$, $y_h$ and $y'_h$. Without lost of generality we can compute $\epsilon$ for  the vectors indexed by $x_e$, $x'_e$, $x_f$, $x'_f$,   $x_g$, $x'_g$, $x_h$, $x'_h$, say, in that order. Hence
\bea
\epsilon=
  \begin{bmatrix}
    1 & 0 & 1 & 0 & 1&0&0&1 \\
    0 & 1 & 0 & 1 & 0&0&0&0 \\
    0 & 0 & 0 & 0 & 0&1&1& 0   
  \end{bmatrix}.
\eea
We deduce the matrix $\tilde\epsilon$ 
\bea\label{eq:matrix}
\tilde\epsilon=
  \begin{bmatrix}
    1 & 0 & 0 \\
    0 & 1 & 0 \\
    0 & 0 & 1   
  \end{bmatrix}.
\eea
From Figure \ref{fig:gluingfpoints} and equation \eqref{eq:matrix}, we have $\tau=3=\rank(\tilde{\epsilon})$. 

The corresponding matrices $\sigma$ and $\tilde{\sigma}$ are given by
\bea
\sigma=
  \begin{bmatrix}
    1 & 0 & 0 & 0 \\
    1 & 0 & 0 & 0 \\
    0 & 0 & 0 & 0 \\
    0 & 0 & 0 & 0 \\
    0 & 0  & 0 & 1 \\ 
    0 & 0 & 0 & 0 \\
    0 & 0 & 0 & 0 \\
    0 & 0 & 0 & 0
  \end{bmatrix},
\eea
and 
\bea
\tilde{\sigma}=
  \begin{bmatrix}
    1 & 0 & 0 & 0 \\
    1 & 0 & 0 & 0 \\
    0 & 0 & 0 & 0 
  \end{bmatrix}.
\eea
\begin{lemma} The number $\tau$ introduced above is
\bea
\tau=\sum_{e\in s}\sum_{f\in \cF_{int}(s)}\sum_{\tilde x_g}\tilde{\epsilon}(f,\tilde x_g)\tilde\sigma(\tilde x_g,e).
\eea
\end{lemma}

\begin{lemma}\label{lemma:partialfincc} Consider a state $\hat{s}$ of the embedded graph $\widehat{\cG}_{\mf^0}$ decomposed into the states $s\in S(\cG_{\mf^0})$ and $s_e\in S^1(H_e)\cup S^2(H_e)$, $e\in\cE$ with the matrices $\tilde{\epsilon}$ and $\tilde\sigma$ as introduced above. Then 
\bea\label{eq:boudaryconnectedd}
 C_\partial(\hat{s})=&&\sum_{e\in\cE}C_\partial(s_e)+C_\partial(s)- |\{s_e\in S^1(H_e)\}|-2|\{s_e\in S^2(H_e)\}|+\cr&&\sum_{e\in s}\sum_{f\in \cF_{int}(s)}\sum_{\tilde x_g}\tilde{\epsilon}(f,\tilde x_g)\tilde\sigma(\tilde x_g,e),
\eea
\bea\label{eq:internalfacee}
F_{\inter}(\hat{s}) =&& \sum_{e\in\cE}F_{\inter}(s_e) + F_{\inter}(s) -  \sum_{e\in s}\sum_{f\in \cF_{int}(s)}\sum_{\tilde x_g}\tilde{\epsilon}(f,\tilde x_g)\tilde\sigma(\tilde x_g,e),
\eea
\bea\label{eq:flafla}
f(\hat{s}) =\sum_{e\in\cE}f(s_e)+f(s)- 2|\{s_e\in S^1(H_e)\}| - 4|\{s_e\in S^2(H_e)\}| .
\eea
\end{lemma}
\proof From the proof of \eqref{eq:boudaryconnectedd} we can deduce \eqref{eq:internalfacee} by using the relation \eqref{rel14} in Lemma \ref{lemma:connectedhats}. Let $e=(m_e,n_e)$, $a_e$ and $a'_e$ be the end points of the arc $m_e$ and $b_e$ and $b'_e$ those of $n_e$. 

- If a connected component of the boundary graph of the graph $s\in S(\cG_{\mf^0})$ is such that it does not contain any of the points $a_e$, $a'_e$, $b_e$ and $b'_e$, then there is a natural corresponding component in $\hat{s}$. 

- If a connected component of the boundary graph of the graph of $s\in S(\cG_{\mf^0})$  contains some of the points $a_e$, $a'_e$, $b_e$ and $b'_e$, then there is a correspondence connected component of the boundary of the graph $\hat{s}\in S(\widehat{\cG}_{\mf^0})$ containing the same set of points. 

- Assume that an internal face of $s\in S(\cG_{\mf^0})$  contains some of the points $a_e$, $a'_e$, $b_e$, $b'_e$. This face can correspond to a connected component of the boundary graph of the graph $\hat{s}\in S(\widehat{\cG}_{\mf^0})$ containing the same set of points. The total number of the components for which some of the points $a_e$, $a'_e$, $b_e$ and $b'_e$ belong to connected component of the boundary graph in $\hat{s}$ is $\sum_{e\in s}\sum_{f\in \cF_{int}(s)}\sum_{\tilde x_g}\tilde{\epsilon}(f,\tilde x_g)\tilde\sigma(\tilde x_g,e)$.

However, $\hat{s}$ has two kind of extra connected components: the unmarked components  in the graphs $s_e$ which do not contain any of the points $a_e$, $a'_e$, $b_e$ and $b'_e$ and the one containing some of these points belonging to connected components of the boundary graphs. The number of the unmarked components  for each $e$ for which $s_e\in  S^1(H_e)$ or $s_e\in  S^1(H_e)$ is  $C_\partial(s_e)-1$ and $C_\partial(s_e)-2$ respectively. 
The number of connected component of the boundary graph of $\hat{s}$ is 
\bea
C_\partial(\hat{s})=&&C_\partial(s) + \sum_{s_e\in S^1(H_e)}\Big(C_\partial(s_e)-1\Big)+ \sum_{s_e\in S^2(H_e)}\Big(C_\partial(s_e)-2\Big)+\cr&&\sum_{e\in s}\sum_{f\in \cF_{int}(s)}\sum_{\tilde x_g}\tilde{\epsilon}(f,\tilde x_g)\tilde\sigma(\tilde x_g,e).
\eea

The proof of \eqref{eq:flafla} is direct since by inserting $s_e\in S^1(H_e)$ in $s$, we keep all the half edges of  $s$ but two half edges of $s_e$ are lost. Inserting  $s_e\in S^2(H_e)$ in $s$ we lose two half edges of $s$ and two half-edges of $s_e$. This ends the proof.
\qed

In order to evaluate the connection between the BR polynomial of $\widehat{\cG}_{\mf^0}$ and those of $\cG_{\mf^0}$ and $H_e$ let us come back to the relations in Lemma \ref{lemma:partialfincc}. It appears that the evaluation of the number of connected components of the boundary graph and internal faces of $\hat{s}$ in term of those in $s_e$ and $s$ becomes more complicated.

We consider the following state sums:
\bea
&&\eta_e^{(1)}(a,b,c,d,l) :=\sum_{s_e\in S^1(H_e)}a^{k(s_e)-1}b^{e(s_e)}c^{F_{\inter}(s_e) }d^{C_\partial(s_e)-1}l^{f(s_e)-2},
\cr
&&\eta_e^{(2)}(a,b,c,d,l) :=\sum_{s_e\in S^2(H_e)}a^{k(s_e)-2}b^{e(s_e)}c^{F_{\inter}(s_e) }d^{C_\partial(s_e)-2}l^{f(s_e)-4}.
\eea
Furthermore let $\mathfrak{F}$ be a map defined for $s\in S(\cG_{f^0})$, by 
\bea\label{eq:fraketa}
\mathfrak{F}\Big(s,\eta_e^{(1)}(a,b,c,d,l)\Big)=\sum_{s_e\in S^1(H_e)}a^{k(s_e)-1}b^{e(s_e)}c^{F_{\inter}(s_e)-\sum_{f\in \cF_{int}(s)}\sum_{\tilde x_g}\tilde{\epsilon}(f,\tilde x_g)\tilde\sigma(\tilde x_g,e) }\times\cr d^{C_\partial(s_e)-1+\sum_{f\in \cF_{int}(s)}\sum_{\tilde x_g}\tilde{\epsilon}(f,\tilde x_g)\tilde\sigma(\tilde x_g,e) }l^{f(s_e)-2},
\eea
where $\epsilon(.)$ and $\tilde\epsilon(.)$ are defined above. 

We remark that if the graph $s$ does not have any internal face i.e $\cF_{int}(s)=\emptyset$ or if  the product $\tilde\epsilon(.)\epsilon(.)$, is always vanishing then the relation \eqref{eq:fraketa} becomes
\bea\label{eq:fraketas}
\mathfrak{F}\Big(s,\eta_e^{(1)}(a,b,c,d,l)\Big)=\eta_e^{(1)}(a,b,c,d,l)\,\,\, \forall s\in S(\cG_{f_0}).
\eea
As consequence, the relation \eqref{eq:fraketa} for  $\eta_e^{(2)}(.)$ gives
\bea\label{eq:fraketad}
\mathfrak{F}\Big(s,\eta_e^{(2)}(a,b,c,d,l)\Big)=\eta_e^{(2)}(a,b,c,d,l)\,\,\, \forall s\in S(\cG_{f_0}).
\eea
Sometimes we will use the following notations:
$$\mathfrak{F}\Big(s,\eta^{(1)}_e(a,b,c,d,l)\Big)=\mathfrak{F}^s_e,\,\, \eta_e^{(1)}(a,b,c,d,l)=\eta_e^{(1)}\,\,\,  \mbox{ and } \,\,\,
\eta_e^{(2)}(a,b,c,d,l)=\eta_e^{(2)}.$$ 

\begin{lemma}\label{lemma:gyvgviyvg}
Let $(\cG_{\mf^0},\{H_e\}_{e\in \cE})$ be a $2-$decomposition of $\widehat{\cG}_{\mf^0}$ and $\mathfrak{F}^s_e$, $\eta_e^{(2)}$ two functions as introduced above. Then
\bea
Z(\widehat{\cG}_{\mf^0};a,b,c,d,l) = &&\sum_{s\in S(\cG_{\mf^0})} a^{k(s)}c^{F_{\inter}(s)}d^{C_\partial(s)}l^{f(s)}\times \cr&&\Big(\prod_{e\in s}\mathfrak{F}^s_e\Big)\Big(\prod_{e\notin s}\eta_e^2\Big).
\eea
\end{lemma}

\proof
We recall that any state $\hat{s}$ is decomposed into states $s\in S(\cG_{\mf^0})$, $s_e\in S^1(H_e)$ and $t_h\in S^2(H_h)$. By Lemma \ref{lemma:partialfincc}, we have 
\bea
&&a^{k(s)}c^{F_{\inter}(s)}d^{C_\partial(s)}l^{f(s)}\prod_{e\in s}a^{k(s_e)-1}b^{e(s_e)}c^{F_{\inter}(s_e) -\sum_{f\in \cF_{int}(s)}\sum_{\tilde x_g}\tilde{\epsilon}(f,\tilde x_g)\tilde\sigma(\tilde x_g,e) }\times\cr&&d^{\partial(s_e)-1+\sum_{f\in \cF_{int}(s)}\sum_{\tilde x_g}\tilde{\epsilon}(f,\tilde x_g)\tilde\sigma(\tilde x_g,e) }l^{f(s_e)-2}\times
\prod_{h\notin s}a^{k(t_h)-1}b^{e(t_h)}c^{F_{\inter}(t_h) }d^{\partial(t_h)-2}l^{f(t_h)-4}
\cr&&=a^{k(s)+\sum (k(s_e)-1)+\sum(k(t_h)-2)} c^{\sum F_{\inter}(s_e)+\sum F_{\inter}(s_h)-\sum_{e\in s}\sum_{f\in \cF_{int}(s)}\sum_{\tilde x_g}\tilde{\epsilon}(f,\tilde x_g)\tilde\sigma(\tilde x_g,e) }\times
\cr&&\quad d^{\sum (C_\partial(s_e)-1)+\sum (C_\partial(t_h)-2)+\sum_{e\in s}\sum_{f\in \cF_{int}(s)}\sum_{\tilde x_g}\tilde{\epsilon}(f,\tilde x_g)\tilde\sigma(\tilde x_g,e)}l^{\sum (f(s_e)-2)+\sum (f(t_h)-4)}\times
\cr&& \quad  b^{\sum e(s_e)+\sum e(t_h)}
\cr&&= a^{k(\hat{s})}b^{e(\hat{s})}c^{F_{\inter}(\hat{s})}d^{C_\partial(\hat{s})}l^{f(\hat{s})}.\nonumber
\eea
\qed

\begin{example}
Let us compute $Z(.)$ for the $2-$decomposition $(G,\{H_f,H_g\})$ of the graph $\hat{G}$ shown in Figure \ref{fig:twosepgraphs}

\begin{figure}[h]
 \centering
     \begin{minipage}[t]{.8\textwidth}
      \centering
\includegraphics[angle=0, width=5.5cm, height=4cm]{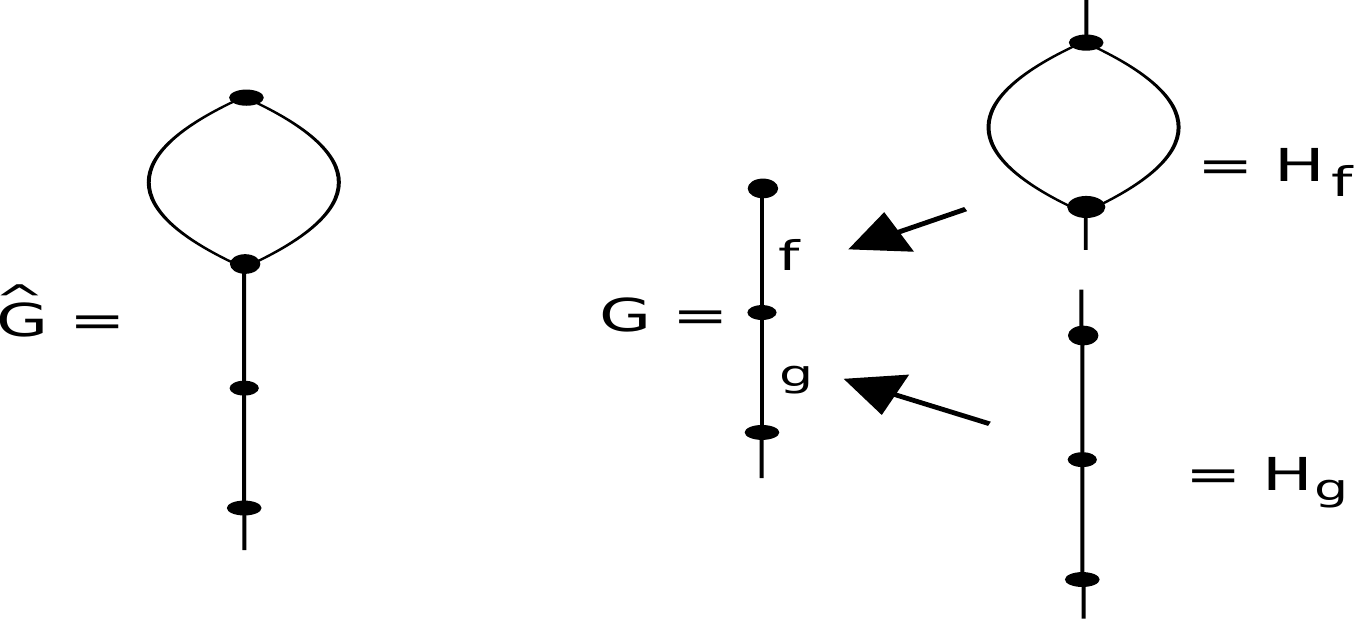}
\caption{ {\small A HERG $\widehat{G}$ on the left and its $2-$decomposition $(G,\{H_f,H_g\})$ on the right}}
\label{fig:twosepgraphs}
\end{minipage}
\end{figure}
 
In this example, $\cF_{\inter}(G)=\emptyset$ and
\bea
Z(\widehat{G};a,b,c,d,l)=adl\eta_f^1\eta_g^1+a^2d^2l^3\eta_f^1\eta_g^2
+a^2d^2l^3\eta_f^2\eta_g^1+a^3d^3l^5\eta_f^2\eta_g^2,
\eea
\bea
\eta_f^1(a,b,c,d,l)=b^2c + 2bl^2,\quad \eta_f^2(a,b,c,d,l)=l^2\cr \eta_g^1(a,b,c,d,l)=b^2,\quad \eta_g^2(a,b,c,d)=2b+adl^2.
\eea
Hence 
\bea
Z(\widehat{G};a,b,c,d,l)&=&adl(b^2c+2bl^2)b^2+a^2d^2l^3(b^2c+2bl^2)(2b+adl^2)+a^2b^2d^2l^2 \cr
&&\quad+a^3d^3l^5(2b+adl^2)\cr
&=&ab^4cdl+2ab^3dl^3+2a^2b^3cd^2l^3+a^3b^2cd^3l^5+4a^2b^2d^2l^5\cr
&&\quad+2a^3bd^3l^7+a^2b^2d^2l^2+2a^3bd^3l^5+a^4d^4l^7.
\eea
\end{example}
\begin{lemma}
Let $(\cG_{\mf^0}, \{H_e\}_{e\in \cE})$ be a two decomposition  of $\widehat{\cG}_{\mf^0}$. We have 
\bea\label{eq:lemprod}
Z(\widehat{\cG}_{\mf^0};a,b,c,d,l)=\Big(\prod_{e\in \cE}\big(\eta_e^{(2)}\big)\Big)Z(\cG_{\mf^0};a,\{\mathfrak{F}^s_e/\eta_e^{(2)}\}_{e\in \cE},c,d,l).
\eea
\end{lemma}
Consider a graph $A_e$ as introduced above with $H_e=A_e\vee e$ and $\widehat{\cG}_{\mf^0}=(\cG_{\mf^0},\{H_e\})$. Let us study the effect of the insertion of the edge $e$ in a state $s$ of $H_e$. If $s\in S^{2}(H_e)$, then the insertion of $e$ in $s$ decreases the number of  connected components and connected components of the boundary graph of $s$ by one and the number of half edges by two. In other words, if $s$ contributes the term $a^{k}b^{e}c^{F_{\inter}}d^{C_\partial}l^f$ in the expression of $Z(.)$, then the state $s\cup e$ contributes the term $a^{k-1}(b^{e}x_e)c^{F_{\inter}}d^{C_\partial-1}l^{f-2}$. Assume that $s\in S^{1}(H_e)$. Three possibilities occur.  In all these possibilities, inserting $e$ in $s\in S^{1}(H_e)$ will decrease by two the number of half edges.

$\bullet$ Inserting $e$ in $s\in S^{1}(H_e)$ increases the number of internal faces by two and decreases by one the number of  connected components of the boundary graph. If $s$ contributes the term $a^{k}b^{e}c^{F_{\inter}}d^{C_\partial}l^f$ to the HERGs BR polynomial then the state $s\cup e$ obtained by inserting $e$ contributes $a^{k}(b^{e}x_e)c^{F_{\inter}+2}d^{C_\partial-1}l^{f-2}$. 

$\bullet$ Inserting $e$ in $s\in S^{1}(H_e)$ increases the number of internal faces by one such that if $s$ contributes the term $a^{k}b^{e}c^{F_{\inter}}d^{C_\partial}l^f$ to the HERGs BR polynomial then the state $s\cup e$ obtained by inserting $e$ contributes $a^{k}(b^{e}x_e)c^{F_{\inter}+1}d^{C_\partial}l^{f-2}$. 

$\bullet$ Inserting $e$ in $s\in S^{1}(H_e)$ increases the number of connected components of the boundary graph by one such that if $s$ contributes the term $a^{k}b^{e}c^{F_{\inter}}d^{C_\partial}l^f$ to the HERGs BR polynomial then the state $s\cup e$ obtained by inserting $e$ contributes $a^{k}(b^{e}x_e)c^{F_{\inter}}d^{C_\partial+1}l^{f-2}$. 

We can summarized all these possibilities by saying that if $s\in S^{1}(H_e)$ and $s$ contributes the term $a^{k}b^{e}c^{F_{\inter}}d^{C_\partial}l^f$ to the HERGs BR polynomial then the state $s\cup e$ obtained by inserting $e$ contributes $a^{k}  (b^{e}x_e)c^{F_{\inter}+\theta(s)}d^{C_\partial+1-\theta(s)}l^{f-2}$; $\theta(s)\in \{0,1,2\}$. 

Based on this we can write $\mathfrak{F}\Big(s,\eta_e^{(1)}\Big)$ as
\bea\label{eq:bbblll}
\mathfrak{F}\Big(s,\eta_e^{(1)}\Big)&=&\mathfrak{F}^0\Big(s,\eta_e^{(1)}\Big)+\mathfrak{F}^1\Big(s,\eta_e^{(1)}\Big)+\mathfrak{F}^2\Big(s,\eta_e^{(1)}\Big),
\eea
with 
\bea
\mathfrak{F}^0\Big(s,\eta_e^{(1)}\Big) &=&\sum_{s\in S^1(H_e)|\theta(s)=0}a^{k(s)-1}b^{e(s)}c^{F_{\inter}(s) -\sum_{f\in \cF_{int}(s)}\sum_{\tilde x_g}\tilde{\epsilon}(f,\tilde x_g)\tilde\sigma(\tilde x_g,e) }\times\cr
&&\quad d^{C_\partial(s)-1+\sum_{f\in \cF_{int}(s)}\sum_{\tilde x_g}\tilde{\epsilon}(f,\tilde x_g)\tilde\sigma(\tilde x_g,e) }l^{f(s)-2},\cr\cr \mathfrak{F}^1\Big(s,\eta_e^{(1)}\Big)&=&\sum_{s\in S^1(H_e)|\theta(s)=1}a^{k(s)-1}b^{e(s)}c^{F_{\inter}(s) -\sum_{f\in \cF_{int}(s)}\sum_{\tilde x_g}\tilde{\epsilon}(f,\tilde x_g)\tilde\sigma(\tilde x_g,e)  }\times\cr &&\quad d^{C_\partial(s)-1+\sum_{f\in \cF_{int}(s)}\sum_{\tilde x_g}\tilde{\epsilon}(f,\tilde x_g)\tilde\sigma(\tilde x_g,e) }l^{f(s)-2},\cr\cr\mathfrak{F}^2\Big(s,\eta_e^{(1)}\Big) &=&\sum_{s\in S^1(H_e)|\theta(s)=2}a^{k(s)-1}b^{e(s)}c^{F_{\inter}(s) -\sum_{f\in \cF_{int}(s)}\sum_{\tilde x_g}\tilde{\epsilon}(f,\tilde x_g)\tilde\sigma(\tilde x_g,e)  }\times\cr &&\quad d^{C_\partial(s)-1+\sum_{f\in \cF_{int}(s)}\sum_{\tilde x_g}\tilde{\epsilon}(f,\tilde x_g)\tilde\sigma(\tilde x_g,e)  }l^{f(s)-2}.
\eea
The map $\mathfrak{F}$ applied on $Z_{A_e}(a,b,c,d,l)$ and $s\in S(\cG_{f^0})$, gives
\bea\label{eq:fcontrade}
\mathfrak{F}\Big(s,Z_{A_e}(a,b,c,d,l)\Big)&=&\sum_{s_e\in S(A_e)}a^{k(s_e)}b^{e(s_e)}c^{F_{\inter}(s_e)-\sum_{f\in \cF_{int}(s)}\sum_{\tilde x_g}\tilde{\epsilon}(f,\tilde x_g)\tilde\sigma(\tilde x_g,e)  }\times 
\cr &&d^{C_\partial(s_e)+\sum_{f\in \cF_{int}(s)}\sum_{\tilde x_g}\tilde{\epsilon}(f,\tilde x_g)\tilde\sigma(\tilde x_g,e)  }l^{f(s_e)}.
\eea
Then we can write the following decomposition
\bea
\mathfrak{F}\Big(s,Z_{A_e}(a,b,c,d,l)\Big)&=&\sum_{s\in S^{1}(H_e)}a^{k(s_e)}b^{e(s_e)}\times \cr &&c^{F_{\inter}(s_e)-\sum_{f\in \cF_{int}(s)}\sum_{\tilde x_g}\tilde{\epsilon}(f,\tilde x_g)\tilde\sigma(\tilde x_g,e)  }\times \cr &&d^{C_\partial(s_e)+\sum_{f\in \cF_{int}(s)}\sum_{\tilde x_g}\tilde{\epsilon}(f,\tilde x_g)\tilde\sigma(\tilde x_g,e) }l^{f(s_e)} \cr &+& x_e\sum_{s\in S^{1}(H_e)|\theta(s_e)=0,1,2}a^{k(s_e)}b^{e(s_e)}\times \cr && c^{F_{\inter}(s_e)-\sum_{f\in \cF_{int}(s)}\sum_{\tilde x_g}\tilde{\epsilon}(f,\tilde x_g)\tilde\sigma(\tilde x_g,e) }\times\cr
&&d^{C_\partial(s_e)+\sum_{f\in \cF_{int}(s)}\sum_{\tilde x_g}\tilde{\epsilon}(f,\tilde x_g)\tilde\sigma(\tilde x_g,e) \,+1-\theta(s_e)}l^{f(s_e)}\cr &+&(1+a^{-1}x_ed^{-1}l^{-2})\sum_{s\in S^{2}(H_e)}a^{k(s_e)}b^{e(s_e)}\times \cr &&c^{F_{\inter}(s_e)-\sum_{f\in \cF_{int}(s)}\sum_{\tilde x_g}\tilde{\epsilon}(f,\tilde x_g)\tilde\sigma(\tilde x_g,e) }\times \cr &&d^{C_\partial(s_e)+\sum_{f\in \cF_{int}(s)}\sum_{\tilde x_g}\tilde{\epsilon}(f,\tilde x_g)\tilde\sigma(\tilde x_g,e) }l^{f(s_e)},
\eea
or
\bea\label{eq:uhg}
\mathfrak{F}\Big(s,Z_{A_e}(a,b,c,d,l)\Big)&=&adl^{2}\mathfrak{F}\big(s,\eta_e^{(1)}\big)+x_el^{2}\Big(ad^2\mathfrak{F}^0\big(s,\eta_e^{(1)}\big) +acd\mathfrak{F}^1\big(s,\eta_e^{(1)}\big) +\mathfrak{F}^2\big(s,\eta_e^{(1)}\big)\Big)\cr
&&\quad+(1+a^{-1}x_ed^{-1}l^{-2})a^2d^2l^{4}\eta_e^{(2)}\cr &=& adl^{2}\Big(\mathfrak{F}\big(s,\eta_e^{(1)}\big)+adl^{2}\eta_e^{(2)}\Big)+ax_el^{2}\Big(d^2\mathfrak{F}^0\big(s,\eta_e^{(1)}\big) +cd\mathfrak{F}^1\big(s,\eta_e^{(1)}\big) \cr
&&+\mathfrak{F}^2\big(s,\eta_e^{(1)}\big)+d\eta_e^{(2)}\Big).
\eea
Furthermore 
\bea
Z_{A_e}(a,b,c,d,l)=Z_{H_e}(a,b,c,d,l)+x_eZ_{A_e/e}(a,b,c,d,l).
\eea

There is a one-to-one correspondence between each state $s'\in S(A_e)$ such that $e\in s'$, and a state $\bar{s'}:=s'/e\in S(A_e/e)$. Let us consider the map
\bea\label{eq:fprimcontrade}
\mathfrak{F}'\Big(s,Z_{A_e/e}(a,b,c,d,l)\Big)=\sum_{s'/e\in S(A_e/e)}a^{k(s'/e)}b^{e(s'/e)}c^{F_{\inter}(s'/e)-\sum_{f\in \cF_{int}(s')}\sum_{\tilde x_g}\tilde{\epsilon}(f,\tilde x_g)\tilde\sigma(\tilde x_g,e)  }\times 
&&\cr d^{C_\partial(s'/e)+\sum_{f\in \cF_{int}(s')}\sum_{\tilde x_g}\tilde{\epsilon}(f,\tilde x_g)\tilde\sigma(\tilde x_g,e)  }l^{f(s'/e)},
\eea
where $\tilde{\epsilon}$ and $\tilde\sigma$ are defined for $s'$.

The following lemma is direct
\begin{lemma}
Given a state $s\in S(\cG_{f^0})$, we have
\bea\label{eq:addi}
\mathfrak{F}(s,Z_{A_e}(a,b,c,d,l)=\mathfrak{F}(s,Z_{H_e}(a,b,c,d,l))+x_e\mathfrak{F}'(s,Z_{A_e/e}(a,b,c,d,l)).
\eea
where $\mathfrak{F}'(s,Z_{A_e/e}(a,b,c,d,l))$ is defined as in \eqref{eq:fprimcontrade}.
\end{lemma}

An identification of \eqref{eq:uhg} and \eqref{eq:addi} gives
\bea\label{eq:compliquedf}
adl^{2}\Big(\mathfrak{F}\big(s,\eta_e^{(1)}\big)+adl^{2}\eta_e^{(2)}\Big)&=&\mathfrak{F}(s,Z_{H_e}(a,b,c,d,l))\cr al^{2}\Big(d^2\mathfrak{F}^0\big(s,\eta_e^{(1)}\big) +cd\mathfrak{F}^1\big(s,\eta_e^{(1)}\big) +\mathfrak{F}^2\big(s,\eta_e^{(1)}\big)+d\eta_e^{(2)}\Big)&=& \mathfrak{F}'(s,Z_{A_e/e}(a,b,c,d,l)).
\eea

A substitution of  \eqref{eq:compliquedf} in \eqref{eq:lemprod} together with \eqref{eq:bbblll} gives
 
\begin{theorem}\label{theo:prod}
Let $\widehat{\cG}_{\mf^0}$ be a half edged embedded graph with $2-$decomposition $(\cG_{\mf^0}, \{H_e\}_{e\in\cE})$, such that each graph $H_e$ is embedded  in a neighbourhood of the edge $e$ of the embedded graph $\cG_{\mf^0}$. Let $A_e$ be a HERG as introduced above. Then

\bea
Z_{\widehat{\cG}_{\mf^0}}(a,b,c,d,l)=(adl^{2})^{-e(\cG_{\mf^0})}\Big(\prod_{e\in \cE} g_e\Big)Z_{\cG_{\mf^0}}(a,\{f_e/g_e\},c,d,l),
\eea

where $f_e$ and $g_e$ are the solutions to 
\bea
adl^{2}g_e+f_e&=&\mathfrak{F}(s,Z_{H_e}(a,b,c,d,l))\cr
adl^2(d-c)\mathfrak{F}^0(s,\eta_e^1)+acl^2(c-d)\mathfrak{F}^2(s,\eta_e^1)+cf_e+g_e&=&\mathfrak{F}'(s,Z_{A_e/e}(a,b,c,d,l)).
\eea
\end{theorem}

This leads to

\begin{corollary}\label{coro:cococ}
Let $\cG_{\mf^0}$ be a HERG, $A$ be a planar HERG \footnote{Planar here means that the underline ribbon graph is planar} and $H=A\vee e.$ Then
\bea\label{coro:}
\cR_{\cG_{\mf^0}\otimes A}(x,y,z,w,t)=h^{n(\cG_{\mf^0})}h'^{\rk(\cG_{\mf^0})}\cR_{\cG_{\mf^0}}(\cG_{\mf^0};\frac{\cR_H(x,y,z,w,t)}{h'},\frac{yh'}{h},z,w,t),
\eea 
where $h$ and $h'$ are the unique solution to
\bea
zwt^2\big(zwt^2(x-1)h+h'\big)&=&\mathfrak{F}(s,\cR_{H}(x,y,z,w,t)),\cr
(x-1)^{-k(H)+1}y^{-v(H)+2}z^{-v(H)+3}t^2(w-z^{-1})\times&&\cr\big(w\mathfrak{F}^0(s,\eta_e^1)-z^{-1}\mathfrak{F}^2(s,\eta_e^1)\big) +zwt^2(yh'+h)&=&\mathfrak{F}'(s,\cR_{A/e}(x,y,z,w,t)).
\eea
\end{corollary}
\proof
From \eqref{eq:boloflaandmult}, we have
\bea
\cR_{\cG_{\mf^0}\otimes A}(x, y, z,w,t) = (x-1)^{-k(\cG_{\mf^0}\otimes A)}(yz)^{-v(\cG_{\mf^0}\otimes A)}Z_{\cG_{\mf^0}\otimes A}((x-1)yz^2, yz, z^{-1},w,t). 
\eea
Applying Theorem \ref{theo:prod}, we obtain
\bea\label{eq:blablabla}
\cR_{\cG_{\mf^0}\otimes A}(x, y, z,w,t) = (x-1)^{-k(\cG_{\mf^0}\otimes A)}(yz)^{-v(\cG_{\mf^0}\otimes A)}g^{e(\cG_{\mf^0})}Z_{\cG_{\mf^0}}((x-1)yz^2,\frac{f}{g}, z^{-1},w,t), 
\eea
where $f$ and $g$ are solutions to
\bea\label{eq:syst}
(x-1)yz^2wt^2\big(f+(x-1)yz^2wt^2g\big)&=&\mathfrak{F}(s,Z_{H}((x-1)yz^2, yz, z^{-1},w,t)),\cr
(x-1)yz^2t^2(w-z^{-1})\big(w\mathfrak{F}^0(s,\eta_e^1)-z^{-1}\mathfrak{F}^2(s,\eta_e^1)\big) &&\cr +(x-1)yz^2wt^2\big(z^{-1}f+g\big)&=&\mathfrak{F}'(s,Z_{A/e}((x-1)yz^2, yz, z^{-1},w,t)).\,\,\,\,\,\,\,\,
\eea
Applying \eqref{eq:boloflaandmult} in \eqref{eq:blablabla}, we obtain
\bea
\cR_{\cG_{\mf^0}\otimes A}(x, y, z,w,t) = (x-1)^{-k(\cG_{\mf^0}\otimes A)}(yz)^{-v(\cG_{\mf^0}\otimes A)}g^{e(\cG_{\mf^0})}\Big(\frac{(x-1)yzg}{f}\Big)^{k(\cG_{\mf^0})}\Big(\frac{f}{g}\Big)^{v(\cG_{\mf^0})}\times \cr\cR_{\cG_{\mf^0}}\Big(\frac{(x-1)yzg+f}{f}, \frac{f}{zg}, z,w,t\Big). 
\eea
Now we set 
\bea
h:=(x-1)^{-k(H)+1}(yz)^{-v(H)+2}g, \quad h':=(x-1)^{-k(H)+1}(yz)^{-v(H)+1}f.
\eea
Using the relations $v(\cG_{\mf^0}\otimes A)=(v(H)-2)e(\cG_{\mf^0})+v(\cG_{\mf^0})$ and $k(\cG_{\mf^0}\otimes A)=(k(H)-1)e(\cG_{\mf^0})+k(\cG_{\mf^0})$, we have
\bea
\cR_{\cG_{\mf^0}\otimes A}(x, y, z,w,t) = h^{n(\cG_{\mf^0})}h'^{\rk(\cG_{\mf^0})}\cR_{\cG_{\mf^0}}\big(\frac{((x-1)h+h'}{h'}, \frac{yh'}{h}, z,w,t\big). 
\eea
Using again \eqref{eq:boloflaandmult}  and the identities $k(H)=k(A/e)$ and $v(A/e)=v(H)-1$, the system \eqref{eq:syst} becomes
 \bea
zwt^2\big(zwt^2(x-1)h+h'\big)&=&\mathfrak{F}(s,\cR_{H}(x,y,z,w,t)),\cr
(x-1)^{-k(H)+1}y^{-v(H)+2}z^{-v(H)+3}t^2(w-z^{-1})\times && \cr\big(w\mathfrak{F}^0(s,\eta_e^1)-z^{-1}\mathfrak{F}^2(s,\eta_e^1)\big)+ zwt^2(yh'+h)&=&\mathfrak{F}'(s,\cR_{A/e}(x,y,z,w,t)).
\eea
\qed

The following corollary gives a way to construct many pair of distinct HERGs with the same BR polynomial.
\begin{corollary}
Let $\cG_{\mf^0}\otimes H$ and $\cG_{\mf^0}'\otimes H$ be two embedded HERGs with the property that each copy of $H$ is embedded in the neighborhood of  an edge. Then if,  $\cR_{\cG_{\mf^0}}(x,y,z,w,t)=\cR_{\cG_{\mf^0}'}(x,y,z,w,t)$,  $\cR_{\cG_{\mf^0}\otimes H}(x,y,z,w,t)=\cR_{\cG_{\mf^0}'\otimes H}(x,y,z,w,t)$
\end{corollary}
\proof
The proof of this corollary is a direct consequence of Corollary \ref{coro:cococ}. Let $z=1$, $w=1$ and $t=1$. Then $T(\cG_{\mf^0})=T(\cG_{\mf^0}')$ where $T$ denote the Tutte polynomial obtained by taking the summation over the spanning cutting subgraphs. Furthermore the rank and nullity can be recovered from this polynomial and $n(\cG_{\mf^0})=n(\cG_{\mf^0}')$, $\rk(\cG_{\mf^0})=\rk(\cG_{\mf^0}')$.
\qed

Coming back to the definition of $\tilde\sigma(.)$ and $\tilde{\epsilon}(.)$,  from now, in order to simplify our results we assume that the template has no internal face or the product  $\tilde\sigma(.)\tilde{\epsilon}(.)$ is always equal to zero.

\subsection{The BR polynomial for HERGs: the general case}
\label{sect:bropenss}

Consider the construction of $\widehat{\cG}_{\mf^0}$ from the $2-$decomposition $(\cG_{\mf^0},\{H_e\}_{e\in\cE})$ at an edge $e=(u_e,w_e)$ of the template $\cG_{\mf^0}$. As discussed earlier, the graph $\widehat{\cG}_{\mf^0}$ is obtained by the identification of the arcs $m_e$ and $n_e$ of the vertices $u_e$ and $w_e$ in $H_e$ with their correspondence in $\cG_{\mf^0}$. It is important to compare again in this subsection the two definitions of $2-$decomposition and discuss a partition of the set $S(H_e)$ in each case.

From Definition \ref{def:graphhat}, the graph $H_e$ is obtained from $A_e$ by deleting the edge $e$, i.e $A_e=H_e-e$. The partition of $S(H_e)$ will be

$\bullet$ $\ddot{S}^2(H_e)$: the set of states of $H_e$ in which $m_e$ and $n_e$ lie in different connected components and different internal faces or connected components of the boundary graph. We may need to split the set $\ddot{S}^2(H_e)$, in two subsets. The first containing the set of states of $H_e$ in which $m_e$ and $n_e$ lie in different connected components and different internal faces. The second containing the set of states of $H_e$ in which $m_e$ and $n_e$ lie in different connected components and different connected components of the boundary graph.

$\bullet$ $\bar{S}^1(H_e)$: the set of states of $H_e$ in which $m_e$ and $n_e$ lie in the same connected components and same internal face or the same connected component of the boundary graph. We may split in the same way as above the set $\bar{S}^1(H_e)$ in two subsets.

$\bullet$ $\ddot{S}^1(H_e)$: the set of states of $H_e$ in which $m_e$ and $n_e$ lie in the same connected components and different connected components of the boundary graph or internal faces. In this case three different possibilities may occur. The points $m_e$ and $n_e$ lie in two different internal faces or two different connected components of the boundary graph. Else one of the points $m_e$ and $n_e$ lies in an internal face and the second in a connected component of the boundary graph.

Finally, Definition \ref{def:graphhat} may give use a partition of $S(H_e)$ in seven different subsets. Now from Definition \ref{def:graphhatclass2}, the partition of $S(H_e)$ will be

$\bullet$ $\ddot{S}^2(H_e)$: the set of states of $H_e$ in which $m_e$ and $n_e$ lie in different connected components and different connected components of the boundary graph.

$\bullet$ $\bar{S}^1(H_e)$: the set of states of $H_e$ in which $m_e$ and $n_e$ lie in the same connected components and same connected components of the boundary graph.

$\bullet$ $\ddot{S}^1(H_e)$: the set of states of $H_e$ in which $m_e$ and $n_e$ lie in the same connected components and different connected components of the boundary graph.

Then we obtain $S(H_e)=\ddot{S}^2(H_e)\cup \bar{S}^1(H_e)\cup \ddot{S}^1(H_e)$. We adopt Definition \ref{def:graphhatclass2} and we have

\begin{lemma}
The sets $\ddot{S}^2(H_e)$, $\bar{S}^1(H_e)$ and $\ddot{S}^1(H_e)$ partition the set $S(H_e)$, $e\in\cE$. Moreover every state of $\widehat{\cG}$ can be uniquely obtained by the replacement of an edge in a state of  $\cG$ by an element of $\bar{S}^1(H_e)$ and a replacement of  a non-edge in that state of $\cG$ in an element of  $\bar{S}^1(H_e)\cup\ddot{S}^1(H_e)$. 
\end{lemma}

If we want to construct the set of states of $S(\widehat{\cG}_{\mf^0})$ by using the states of $\cG_{\mf^0}$, we run into a problem. We then introduce a class of  states $\widetilde{\cG}_{\mf^0}=\cG_{\mf^0}\otimes T$ where $T=\includegraphics[angle=0, width=1.2cm, height=1.5cm]{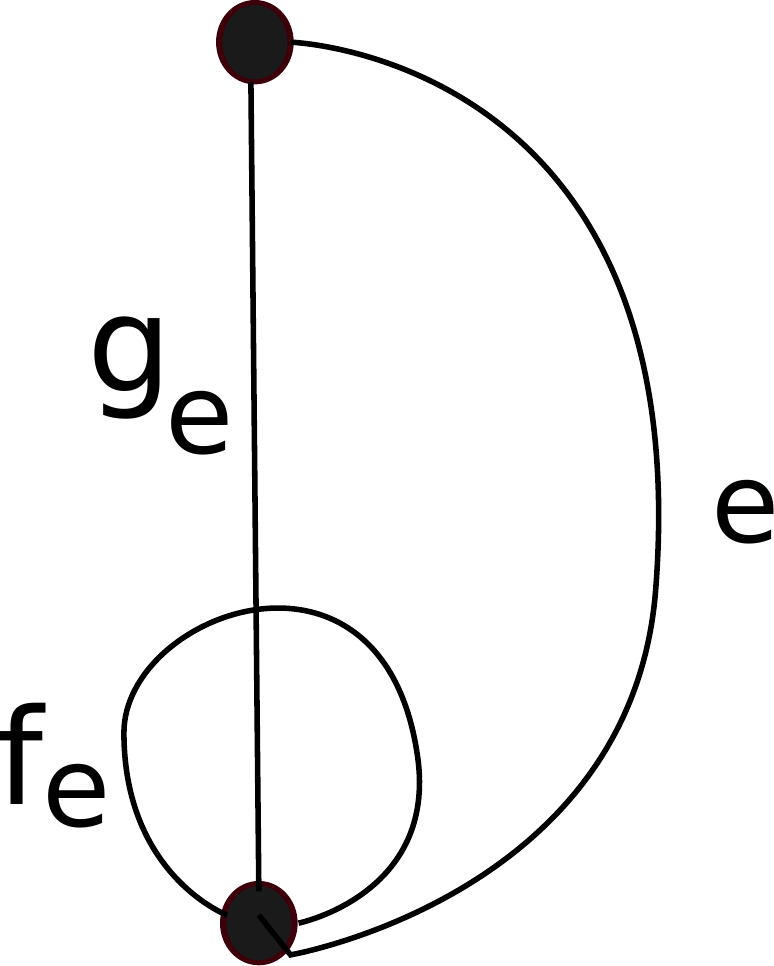}$. We say that two states of $\widetilde{\cG}_{\mf^0}$ are equivalent if for some choice of $e$, one state contains neither of the edges $g_e$ or $f_e$, the other state contains the edge $f_e$ but not $g_e$, and the remaining edges contained in the two states are the same. This gives an equivalence relation $\sim$ and we say that two states are equivalent if they contain both of $f_e$ and $g_e$ or they contain both $g_e$ but not $f_e$ or one contains neither $g_e$ or $f_e$  and the second does not contain $g_e$ but may or may not contain $f_e$.

Consider the set $S(\widetilde{\cG}_{\mf^0})/\sim$ and $[s]\in S(\widetilde{\cG}_{\mf^0})/\sim$ a class with $s$ one of its representative. A state of $\widehat{\cG}_{\mf^0}$ is obtained by the following 

$\bullet$ if the state $s$ contains both of $f_e$ and $g_e$, we remove both of them and replace in a state from $\ddot{S}^1(H_e)$. 

$\bullet$ if the state $s$ contains $g_e$ but not $f_e$, we remove $g_e$ and replace in a state from $\bar{S}^1(H_e)$. 

$\bullet$ if the state $s$ contains $f_e$ but not $g_e$, we remove $f_e$ and replace in a state from $\ddot{S}^2(H_e)$, or if $s$ contains neither of the edges $g_e$ or $f_e$, then we replace in a state from $\ddot{S}^2(H_e)$.

\begin{lemma}
The sets $\ddot{S}^2(H_e)$, $\bar{S}^1(H_e)$ and $\ddot{S}^1(H_e)$ partition the set $S(H_e)$, $e\in\cE$. Moreover every state of $\widehat{\cG}_{\mf^0}$ can be uniquely obtained by replacing classes in  $S(\widetilde{\cG}_{\mf^0})/\sim$ by an element of $\ddot{S}^2(H_e)$, $\bar{S}^1(H_e)$ and $\ddot{S}^1(H_e)$ in the manner described above.
\end{lemma}

Let us consider the following sums:
\bea
&&\Phi_{\widetilde{\cG}}(a,\{g_e,f_e\}_{e\in\cE},c,d) :=\sum_{[s]\in S(\widetilde{\cG}_{\mf^0})/\sim}a^{k(s)}\big(\prod_{e\in s}x_e\big)c^{F_{\inter}(s) }d^{C_\partial(s)},
\cr
&&\ddot{\eta}_e^{(1)}(a,b,c,d) :=\sum_{s\in \ddot{S}^1(H_e)}a^{k(s)-1}b^{e(s)}c^{F_{\inter}(s) }d^{C_\partial(s)-2},
\cr
&&\ddot{\eta}_e^{(2)}(a,b,c,d) :=\sum_{s\in \ddot{S}^2(H_e)}a^{k(s)-2}b^{e(s)}c^{F_{\inter}(s) }d^{C_\partial(s)-2},
\cr
&&\bar{\eta}_e^{(1)}(a,b,c,d) :=\sum_{s\in \bar{S}^1(H_e)}a^{k(s)-1}b^{e(s)}c^{F_{\inter}(s) }d^{C_\partial(s)-1},
\eea
where the product in $\Phi_{\widetilde{\cG}_{\mf^0}}(.)$ is over a representative $s$ of $[s]$. For simplicity we suppose that $s$ has the fewest edges in its equivalence class and $x_e$ stand for $(f_e,g_e)$.

\begin{lemma}\label{lemma:partialfinccnsd}
Suppose that a state $\hat{s}$ of $\widehat{\cG}_{\mf^0}$ is obtained by the replacement from the states $[s]\in S(\widetilde{\cG}_{\mf^0})$, $s_e\in \ddot{S}^1(H_e)$, $t_e\in \ddot{S}^1(H_e)$, and $u_e\in \bar{S}^1(H_e)$ using the decomposition above. Then
\bea\label{eq:blablasw}
&&k(\hat{s})=\sum (k(s_e) - 1)+\sum (k(t_e) - 1) \sum (k(u_e) - 1)+ k(s),
\cr
&& C_\partial(\hat{s})=\sum (C_\partial(s_e) - 2)+ \sum (C_\partial(t_e) - 2)+\sum (C_\partial(u_e) - 1)+C_\partial(s),
\cr
&&F_{\inter}(\hat{s}) =\sum F_{\inter}(s_e) + \sum F_{\inter}(t_e) + \sum F_{\inter}(u_e)+F_{\inter}(s),
\cr
&& e(\hat{s})=\sum e(s_e)+\sum e(t_e)+\sum e(u_e).
\eea
where the representative $s\in[s]$ has a fewest number of edges in its class.
\end{lemma}

\proof
The proof of this lemma will follow the one given in Lemma \ref{lemma:gyvgviyvg}. Consider a representative $s\in[s]\in S(\widetilde{\cG}_{\mf^0})/\sim$ which has the fewest number of edges. A connected component of $s$ containing any of the points $a_e$, $a'_e$, $b_e$ and $b'_e$ corresponds to a connected component of the boundary of the graph $\hat{s}$ containing the same set of points. The extra connected components of the boundary graph of $\hat{s}$ arise from the states $s_e\in \ddot{S}^1(H_e)$, $t_e\in \ddot{S}^1(H_e)$, and $u_e\in \bar{S}^1(H_e)$. These components are precisely the unmarked components. For each $e$ they are  $(C_\partial(s_e) - 2)$, $(C_\partial(t_e) - 2)$ and $(C_\partial(u_e) - 1)$ extra components. This ends the proof of the second relation in \eqref{eq:blablasw}. The remaining follow.
\qed
 
Let us consider the map $\mathcal{F}:$ $\mathbb{Z}[\{f_e,g_e\}_{e\in\cE}]$ $\rightarrow$ $\mathbb{Z}[\{\ddot{\eta}_e^{(1)},\ddot{\eta}_e^{(1)},\bar{\eta}_e^{(1)}\}_{e\in\cE}]$ as the linear extension of 

\bea
\mathcal{F}: \prod_{e\in \cE}f_e^{\alpha_e}g_e^{\beta_e} \mapsto \prod_{e\in \cE}(\ddot{\eta}_e^{(1)})^{\alpha_e\beta_e}(\ddot{\eta}_e^{(2)})^{1-\beta_e}(\bar{\eta}_e^{(1)})^{\beta_e-\alpha_e\beta_e}.
\eea

\begin{lemma}
The multivariate version of BR polynomial on HERGs, $Z$, of the graph $\widehat{\cG}_{\mf^0}$ introduced above is given by
\bea
Z_{\widehat{\cG}_{\mf^0}}(a,b,c,d)=\mathcal{F}(\Phi_{\widetilde{\cG}_{\mf^0}}).
\eea
\end{lemma}

\proof
Each state $\hat{s}$ of $\widehat{\cG}_{\mf^0}$ is obtained in the following way: we consider a representative $s\in[s]\in S(\widetilde{\cG}_{\mf^0})/\sim$ which has the fewest number of edges. We replace the three edges configurations $\{f_e,g_e\}$, $g_e$ and $\emptyset$ of the pair of edges $\{f_e,g_e\}$ (as introduced above) by $s_e\in \ddot{S}^1(H_e)$, $\ddot{S}^1(H_e)$, and $\bar{S}^1(H_e)$ respectively. The contribution of $\hat{s}$ in $Z_{\widehat{\cG}_{\mf^0}}(.)$ is 
\bea\label{eq:prodn}
a^{k(s)}c^{F_{\inter}(s) }d^{C_\partial(s)}\prod_{e\in \cE}\big(a^{k(s_e)-1}b^{e(s_e)}c^{F_{\inter}(s_e) }d^{C_\partial(s_e)-2}\big)^{\alpha_e\beta_e}\big(a^{k(s_e)-2}b^{e(s_e)}c^{F_{\inter}(s_e) }d^{C_\partial(s_e)-2}\big)^{1-\beta_e}\cr\times\big(a^{k(s_e)-1}b^{e(s_e)}c^{F_{\inter}(s_e) }d^{C_\partial(s_e)-1}\big)^{\beta_e-\alpha_e\beta_e}=\mathcal{F}\big(a^{k(s)}\big(\prod_{e\in s}x_e\big)c^{F_{\inter}(s) }d^{C_\partial(s)}\big),
\eea
where $\alpha_e = \left\{\begin{array}{ll} 1 &  {\text{if}}\,\, f_e\in s,\\0 & {\text{otherwise}},
\end{array} \right.$ and $\beta_e = \left\{\begin{array}{ll} 1 &  {\text{if}}\,\, g_e\in s,\\0 & {\text{otherwise}}\end{array} \right.$. The expression \eqref{eq:prodn} is also equal to
\bea\label{eq:expression}
a^{k(s)+\sum (k(s_e) - 1)+\sum (k(t_e) - 1) +\sum (k(u_e) - 1)}b^{\sum e(s_e)+\sum e(t_e)+\sum e(u_e)}\times \cr c^{\sum F_{\inter}(s_e) + \sum F_{\inter}(t_e) + \sum F_{\inter}(u_e)+F_{\inter}(s)} d^{\sum (C_\partial(s_e) - 2)+ \sum (C_\partial(t_e) - 2)+\sum (C_\partial(u_e) - 1)+C_\partial(s)},
\eea
where $s_e\in \ddot{S}^1(H_e)$, $t_e\in \ddot{S}^1(H_e)$, and $u_e\in \bar{S}^1(H_e)$. Using Lemma \ref{lemma:partialfinccnsd}, we deduce that the relation in \eqref{eq:expression} is equal to $a^{k(\hat{s})}b^{e(\hat{s})}c^{F_{\inter}(\hat{s})}d^{C_\partial(\hat{s})}$. Since the decomposition of $\hat{s}$ into $s$ and $s_e$ is unique, the result follows summing over the states.
\qed

Remark that the proof of this lemma follows the proof of Lemma 8 in \cite{stephen-fathom}. Following the proof of Lemma 9 in the same paper, we can prove the following result.

\begin{lemma}
 Let $(\cG_{\mf^0},\{H_e\}_{e\in\cE})$ be a two decomposition of $\hat{\cG_{\mf^0}}$, and $\mathcal{H}:$ $\mathbb{Z}[\{f_e,g_e\}_{e\in\cE}]$ $\rightarrow$ $\mathbb{Z}[\{\ddot{\eta}_e^{(1)},\ddot{\eta}_e^{(1)},\bar{\eta}_e^{(1)}\}_{e\in\cE}]$ be the linear extension of  the map
\bea
\mathcal{H}: \prod_{e\in \cE}f_e^{\alpha_e}g_e^{\beta_e} \mapsto \prod_{e\in \cE}(\ddot{\eta}_e^{(1)})^{\alpha_e\beta_e}(\ddot{\eta}_e^{(2)})^{1-\beta_e}(\bar{\eta}_e^{(1)})^{\beta_e-\alpha_e\beta_e}.\nonumber
\eea
Then
\bea
Z(\widehat{\cG}_{\mf^0};a,b,c,d)=\mathcal{H}(Z(\widetilde{\cG}_{\mf^0};a,x,c,d)).
\eea
\end{lemma}

Let us study the effect of the insertion of the edge $e$ in a state of $H_e$: $s\in \ddot{S}^1(H_e)\cup \ddot{S}^2(H_e)\cup \bar{S}^1(H_e)$. The effect of inserting $s$ in $\bar{S}^1(H_e)\cup \ddot{S}^2(H_e)$ is similar to the one we discussed in the planar case. If $s\in \ddot{S}^{2}(H_e)$, then the insertion of $e$ in $s$ decreases the number of  connected components and connected components of the boundary graph by one. This means, if $s$ contributes the term $a^{k}b^{e}c^{F_{\inter}}d^{C_\partial}$, then the state $s\cup e$ contributes the term $a^{k-1}(b^{e}x_e)c^{F_{\inter}}d^{C_\partial-1}$. The case $s\in \bar{S}^1(H_e)$ is summarized in three possibilities: if $s$ contributes the term $a^{k}b^{e}c^{F_{\inter}}d^{C_\partial}$ to the HERGs BR polynomial then the state $s\cup e$ obtained by inserting $e$ contributes $a^{k} (b^{e}x_e)c^{F_{\inter}+\theta(s)}d^{C_\partial+1-\theta(s)}$; $\theta(s)\in \{0,1,2\}$.  Assume that $s\in \ddot{S}^{1}(H_e)$. Tow possibilities occur.  

$\bullet$ Inserting $e$ in $s\in \ddot{S}^{1}(H_e)$ increases the number of internal faces by one and decreases by two the number of  connected components of the boundary graph. If $s$ contributes the term $a^{k}b^{e}c^{F_{\inter}}d^{C_\partial}$ to the HERGs BR polynomial then the state $s\cup e$ obtained by inserting $e$ contributes $a^{k}(b^{e}x_e)c^{F_{\inter}+1}d^{C_\partial+2}$. 

$\bullet$ Inserting $e$ in $s\in \ddot{S}^{1}(H_e)$ decreases the number of connected components of the boundary graph by one such that if $s$ contributes the term $a^{k}b^{e}c^{F_{\inter}}d^{C_\partial}$ to the HERGs BR polynomial then the state $s\cup e$ obtained by inserting $e$ contributes $a^{k}(b^{e}x_e)c^{F_{\inter}}d^{C_\partial-1}$. 

We can summarized our discussion by saying that is $s\in \ddot{S}^{1}(H_e)$ contributes the term $a^{k}b^{e}c^{F_{\inter}}d^{C_\partial}$ to the HERGs BR polynomial, then the state $s\cup e$ obtained by inserting $e$ in $s$ contributes 
$a^{k}(b^{e}x_e)c^{F_{\inter}+\theta(s)}d^{C_\partial-1-\theta(s)}$; $\theta(s)\in \{0,1\}$.

Based on the previous discussions, we can write $\bar{\eta}_e^{(1)}(a,b,c,d)$ as
\bea
\bar{\eta}_e^{(1)}(a,b,c,d) &=&\bar{\eta}_e^{(1,0)}(a,b,c,d) +\bar{\eta}_e^{(1,1)}(a,b,c,d) +\bar{\eta}_e^{(1,2)}(a,b,c,d),\cr \ddot{\eta}_e^{(1)}(a,b,c,d) &=&\ddot{\eta}_e^{(1,0)}(a,b,c,d) +\ddot{\eta}_e^{(1,1)}(a,b,c,d),
\eea
with 
\bea
&&\bar{\eta}_e^{(1,0)}(a,b,c,d) =\sum_{s\in \bar{S}^1(H_e)|\theta(s)=0}a^{k(s)-1}b^{e(s)}c^{F_{\inter}(s) }d^{C_\partial(s)-1},\cr &&\bar{\eta}_e^{(1,1)}(a,b,c,d) =\sum_{s\in \bar{S}^1(H_e)|\theta(s)=1}a^{k(s)-1}b^{e(s)}c^{F_{\inter}(s) }d^{C_\partial(s)-1},\cr &&\bar{\eta}_e^{(1,2)}(a,b,c,d) =\sum_{s\in \bar{S}^1(H_e)|\theta(s)=2}a^{k(s)-1}b^{e(s)}c^{F_{\inter}(s) }d^{C_\partial(s)-1},\cr &&\ddot{\eta}_e^{(1,0)}(a,b,c,d) =\sum_{s\in \ddot{S}^1(H_e)|\theta(s)=0}a^{k(s)-1}b^{e(s)}c^{F_{\inter}(s) }d^{C_\partial(s)-2},\cr &&\ddot{\eta}_e^{(1,1)}(a,b,c,d) =\sum_{s\in \bar{S}^1(H_e)|\theta(s)=1}a^{k(s)-1}b^{e(s)}c^{F_{\inter}(s) }d^{C_\partial(s)-2}.
\eea
We also have
\bea
Z_{A_e}(a,b,c,d)&=&\sum_{s\in \bar{S}^{1}(H_e)}a^{k(s)}b^{e(s)}c^{F_{\inter}(s)}d^{C_\partial(s)}\cr &+& x_e\sum_{s\in \bar{S}^{1}(H_e)|\theta(s)=0,1,2}a^{k(s)}b^{e(s)}c^{F_{\inter}(s)+\theta(s)}d^{C_\partial(s)+1-\theta(s)} \cr &+&\sum_{s\in \ddot{S}^{1}(H_e)}a^{k(s)}b^{e(s)}c^{F_{\inter}(s)}d^{C_\partial(s)}\cr &+& x_e\sum_{s\in \ddot{S}^{1}(H_e)|\theta(s)=0,1}a^{k(s)}b^{e(s)} c^{F_{\inter}(s)+\theta(s)}d^{C_\partial(s)-1-\theta(s)} \cr &+& (1+a^{-1}x_ed^{-1})\sum_{s\in \ddot{S}^{2}(H_e)}a^{k(s)}b^{e(s)}c^{F_{\inter}(s)}d^{C_\partial(s)}.
\eea

\bea
Z_{A_e}(a,b,c,d)&=&ad\bar{\eta}_e^{(1)}+x_e\big(ad^2\bar{\eta}_e^{(1,0)} +acd\bar{\eta}_e^{(1,1)} +ac^2\bar{\eta}_e^{(1,2)}\big)+ad^2\ddot{\eta}_e^{(1)}\cr&+&x_e\big(ad\ddot{\eta}_e^{(1,0)} +acd^2\ddot{\eta}_e^{(1,1)}\big)+(1+a^{-1}x_ed^{-1})a^2d^2\ddot{\eta}_e^{(2)}\cr &=& ad(\bar{\eta}_e^{(1)}+d\ddot{\eta}_e^{(1)}+ad\ddot{\eta}_e^{(2)})\cr&+&ax_e\big(d^2\bar{\eta}_e^{(1,0)} +cd\bar{\eta}_e^{(1,1)} +c^2\bar{\eta}_e^{(1,2)}+d\ddot{\eta}_e^{(1,0)} +cd^2\ddot{\eta}_e^{(1,1)} +d\ddot{\eta}_e^{(2)}\big).
\eea

Furthermore 
\bea
Z_{A_e}(a,b,c,d)=Z_{H_e}(a,b,c,d)+x_eZ_{A_e/e}(a,b,c,d).
\eea
Then

\bea\label{eq:complique}
ad(\bar{\eta}_e^{(1)}+d\ddot{\eta}_e^{(1)}+ad\ddot{\eta}_e^{(2)})&=&Z_{H_e}(a,b,c,d),\cr a\big(d^2\bar{\eta}_e^{(1,0)} +cd\bar{\eta}_e^{(1,1)} +c^2\bar{\eta}_e^{(1,2)}+d\ddot{\eta}_e^{(1,0)} +cd^2\ddot{\eta}_e^{(1,1)} +d\ddot{\eta}_e^{(2)}\big)&=& Z_{A_e/e}(a,b,c,d).
\eea

The second relation in \eqref{eq:complique} can also be written as
\bea
ad(d-c)\bar{\eta}_e^{(1,0)} + ac(c-d)\bar{\eta}_e^{(1,2)}  + ad(1-cd)\ddot{\eta}_e^{(1,0)} +ad( c\bar{\eta}_e^{1}+cd\ddot{\eta}_e^{1}+\ddot{\eta}_e^{2})= Z_{A_e/e}(a,b,c,d).
\eea
\begin{theorem}
\label{theo:import}
Let $\widehat{\cG}_{\mf^0}$ be a ribbon graph with the $2-$decomposition $(\cG_{\mf^0},\{H_e\})$, and let $A_e$ be the ribbon graph $H_e$ with an additional ribbon $e$ joining the vertices $u_e$ and $w_e$. Then
\bea
Z(\widehat{\cG}_{\mf^0};a,x,c,d)=\cH(Z(\widetilde{\cG}_{\mf^0};a,x,c,d)),
\eea
where $p_e$, $q_e$ and $r_e$ are uniquely determine by the pair of equations
\bea
p_e+adq_e+dr_e&=&Z_{H_e}(a,b,c,d)\cr ad(d-c)\bar{\eta}_e^{(1,0)} + ac(c-d)\bar{\eta}_e^{(1,2)}  + ad(1-cd)\ddot{\eta}_e^{(1,0)} + cp_e+q_e+cdr_e&=& Z_{A_e/e}(a,b,c,d),
\eea

where the $\eta$ are introduced above and $\cH$ is deduced by
\bea
\mathcal{H}: \prod_{e\in \cE}f_e^{\alpha_e}g_e^{\beta_e} \mapsto \prod_{e\in \cE}\big(\frac{r_e}{ac}\big)^{\alpha_e\beta_e}\big(\frac{q_e}{2ac}\big)^{1-\beta_e}c^{\alpha_e\beta_e-\alpha_e}\big(\frac{p_e}{ac})^{\beta_e-\alpha_e\beta_e}.
\eea
\end{theorem}

\section{The rank $D-$weakly colored stranded graph $2$-decomposition and polynomial invariant}  \label{sect:weaklycoloreddec}
\subsection{Weakly colored stranded graph}

In this subsection we assume that the reader is familiar with stranded graphs.
We briefly review here the weakly-colored stranded graphs introduced in \cite{Avo16}.

 A graph $\cG(\cV,\cE)$ is stranded when its vertices and edges are stranded. A rank $D$ stranded vertex is a trivial disc or a  chord diagram with a collection of $2n$ points on the unit circle (called the vertex frontier) paired by $n$ chords and drawn in a specific way. A rank $D$ stranded edge is a collection of segments called  strands respecting some conditions as introduced in \cite{Avo16}.

A rank $D$ stranded graph $\cG$ is a graph $\cG(\cV,\cE)$ which admits: rank $D$ stranded vertices, rank at most $D$ stranded edges such that the intersections of vertices and edges are pairwise distinct. The graph $\cG$ is a rank $D$ tensor graph if the vertices of $\cG$ have a fixed coordination $D+1$ and their  pre-edges have a fixed cardinal $D$. The merged point graph is $K_{D+1}$ and the edges of $\cG$ are of rank $D$. 

Collapsing the stranded vertices of $\cG$ to points and edges to simple lines, the resulting object is a graph.  The graph $\cG$ is said to be connected if its corresponding collapsed graph is connected. From this point, stranded vertices and edges are always connected. If we assign a color from  the set $\{0,\cdots,D\}$ to each edge of $\cG$ such that no two adjacent edges share the same color, then the graph $\cG$ is called a $(D+1)$ colored graph. It is a bipartite graph if the set $\cV$ of vertices is split into two disjoint sets, i.e. $\cV = \cV^+ \cup \cV^-$ with $\cV^+\cap\cV^-= \emptyset$,  such that each edge connects  a vertex $v^+ \in \cV^+$ and a vertex $v^- \in \cV^-$. The rank $D\geq 1$ tensor  graph $\cG$ is  a rank $D\geq 1$-Colored  tensor graph if it is a $(D+1)$-colored and bipartite graph.

The collapsed graph associated to a colored tensor graph is obtained by regarding the tensor graph as a simple bipartite colored graph and is called {\emph{compact}} in the following.
In such a colored graph, there are some informations that we address now.
\begin{definition}[$p$-bubbles \cite{Gurau:2009tz}]
\label{def:pbub}
Let $\cG$ be a rank $D$ colored  tensor graph. 

- A 0-bubble is a vertex of $\cG$

- A 1-bubble is an edge of $\cG$

- For all $p \geq 2$, a $p$-bubble of $\cG$ with colors $i_1 < i_2 < \dots<i_p $, $p\leq D$,
and  $i_k\in \{0, \dots, D\}$ is a connected rank $p-1$ colored tensor graph the compact form of which is a connected  subgraph 
of the compact form of $\cG$ made of edges of colors $\{i_1, \dots, i_p\}$. 
 \end{definition}

The concept of half-edge find a natural extension on colored tensor graphs and allows to  recall the operation called ``cut'' of an edge.

\begin{figure}[h]
 \centering
     \begin{minipage}[t]{.8\textwidth}
      \centering
\includegraphics[angle=0, width=5cm, height=1cm]{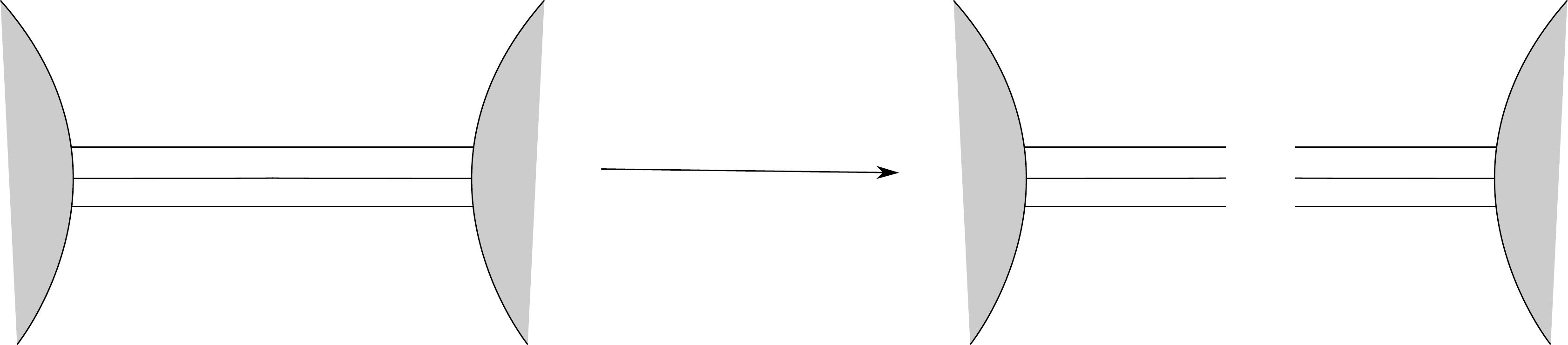}
\caption{ {\small The cut of a rank 3 stranded edge.}}
\label{fig:cuttens}
\end{minipage}
\end{figure}

In the following, a stranded graph having stranded half-edges is denoted by $\cG_{\mf^0}(\cV,\cE)$ or simply $\cG_{\mf^0}$, where $\mf^0$ is the set of the half-edges. 
\begin{definition}[Cut of an edge \cite{riv}]\label{def:cutedtens}
Let $\cG_{\mf^0}(\cV,\cE)$ be a rank D stranded graph
and $e$ a rank $d$ edge of  $\cG_{\mf^0}$, $1\leq d \leq D$.
The cut graph $\cG_{\mf^0} \vee e$ or the graph obtained from $\cG_{\mf^0}$ by cutting $e$ is obtained by replacing the edge $e$ by two rank $d$ stranded half-edges at the end vertices of $e$ and respecting the strand structure
of $e$.  (See Figure \ref{fig:cuttens}.) If $e$ is a loop, the two stranded half-edges are on the same vertex.  
\end{definition}

Using the operation of ``cutting'' of an edge, we obtain a c-subgraph $A_{\mf^0_A}(\cV_A,\cE_A)$ of a rank $D$ stranded graph $\cG_{\mf^0}(\cV,\cE)$ by cutting a subset of edges of $\cG_{\mf^0}$. A spanning c-subgraph called also state $s$ of $\cG_{\mf^0}$ is defined as  a c-subgraph $s_{\mf^0_s}(\cV_s,\cE_s)$ of $\cG_{\mf^0}$ with all vertices  and all additional half-edges of $\cG_{\mf^0}$. Then $\cE_s\subseteq \cE$ and $\cV_s = \cV$, $\mf^0_s= \mf^{0} \cup \mf^{0;1}_s(\cE_s)$, where $\mf^{0;1}_s (\cE_s)$ is the set of half-edges obtained by cutting all edges in $\cE_s'$ (the set of edges incident to the vertices of $s$
and not contained in $\cE_s$) and incident to the vertices of $s$. We write $s \sset  \cG_{\mf^0}$ to say that $s$ is a spanning c-subgraph of $\cG_{\mf^0}$.

The cutting of an edge modifies the strand structure of the graph. In fact as discussed earlier, the presence of half-edges immediately introduce another type of faces which pass through the external points of the half-edges. We then distinguish two kind of faces: open faces which are passing through the external points of the half-edges and the others called closed or internal faces. We denote by $\cF_{\inter}$, the set of closed faces and $\cF_{\ext}$  the set of  open faces.  The set of faces $\cF$ for a rank $D$ half-edged colored tensor graph is then the disjoint union $\cF_{\inter} \cup \cF_{\ext}$. A bubble is called open or external if it contains open faces otherwise it is  closed or internal. We denote respectively by $\cB_{\inter}$ and $\cB_{\ext}$ the sets of closed and open bubbles for a rank $D$ tensor graph.

Their is a graph directly associated to a color tensor graph called  ``boundary graph'' which is obtained by  setting a vertex to each half-edge \cite{Gurau:2009tz}. The boundary graph $\bG(\bV,\bE)$ of a rank $D$ half-edged colored tensor graph  $\cG_{\mf^0}(\cV,\cE)$ is obtained by inserting a vertex with degree $D$ at each additional stranded half-edge of  $\cG_{\mf^0}$ and taking the external faces of $\cG_{\mf^0}$ as its edges. Thus, $|\bV| = |\mf^0|$ and $\bE = \cF_{\ext}$. If the  rank $D$ half-edged colored tensor graph is closed, then its boundary is empty. 

\begin{definition}[Equivalence class of half-edged stranded graph \cite{Avo16}]
Let $D_{\cG_{\mf^0}}$ be the subgraph in a rank $D$ half-edge stranded graph   $\cG_{\mf^0}$ defined by all of its trivial disc vertices and  $\cG_{\mf^0} \setminus D_{\cG_{\mf^0}}$ the rank $D$ half-edges stranded graph obtained after removing  $D_{\cG_{\mf^0}}$ from  $\cG_{\mf^0}$. 

Two rank $D$ half-edged stranded graphs  $\cG_{1,\mf^0(\cG_1)}$ and  $\cG_{2,\mf^0(\cG_2)}$ are ``equivalent up to trivial discs''  if and only if 
$ \cG_{1,\mf^0(\cG_1)}\setminus D_{ \cG_{1,\mf^0(\cG_1)}} = { \cG_{2,\mf^0(\cG_2)}}\setminus D_{ \cG_{2,\mf^0(\cG_2)}}$. 
We note   $ \cG_{1,\mf^0(\cG_1)} \sim { \cG_{2,\mf^0(\cG_2)}}$. 
\end{definition}

As  a consequence of this definition, the contraction of all edges in arbitrary order of a half-edged tensor graph $\cG_{\mf^0}$ yields a half-edged stranded graph  $\cG^0_{\mf^0}$ determined by the boundary $\partial( \cG_{\mf^0})$ up to additional discs. Noting that contracting an edge in a rank $D$ (colored)  half-edged tensor graph does not change its boundary. 

We can now address a precise definition of a rank $D$ w-colored graph.

\begin{figure}[h]
 \centering
     \begin{minipage}[t]{.8\textwidth}
      \centering
\includegraphics[angle=0, width=7cm, height=2cm]{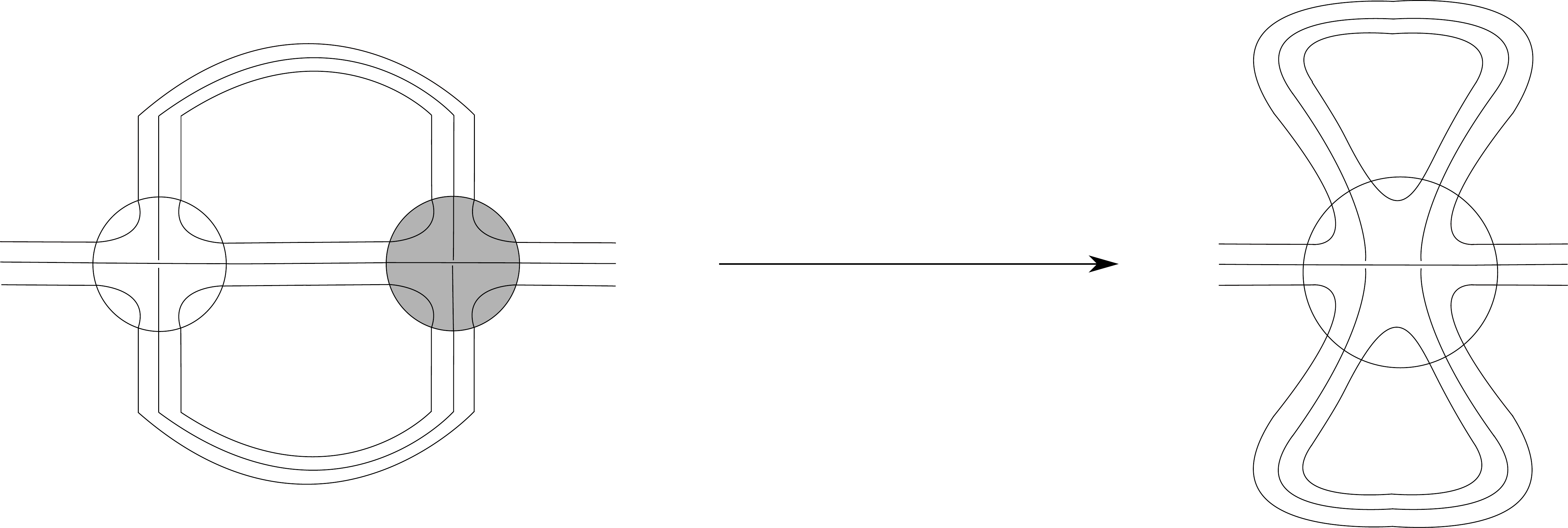}
\vspace{0.1cm}
\caption{ {\small Contraction of an edge in a rank 3 colored tensor graph (left) giving a rank 3 w-colored tensor graph (right)}}
\label{fig:contens}
\end{minipage}
\put(-188,25){1}
\put(-280,25){1}
\put(-233,-5){2}
\put(-233,33){3}
\put(-233,55){0}
\put(-125,25){1}
\put(-68,25){1}
\put(-95,-10){2}
\put(-95,60){0}
\end{figure}
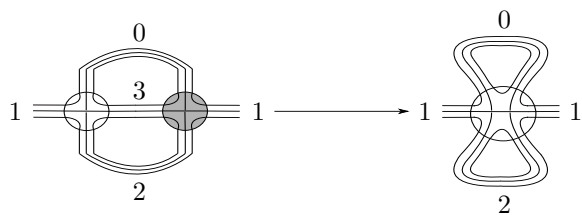

\begin{definition}[Rank $D$ w-colored graph \cite{Avo16}]\label{wdcolo}
A rank $D$ weakly colored  or w-colored graph is the equivalence class (up to trivial discs) of a rank $D$ half-edged stranded graph obtained by successive edge contractions of some rank $D$ half-edged colored tensor graph. An  illustration is given in Figure \ref{fig:contens} .
\end{definition}

\subsection{Polynomial invariant expansion}

Here again, some definitions and notations introduced in \cite{Avo16} deserve to be review as well. 

Consider a representative $ \cG_{\mf^0}$ of any rank $D$ w-colored graph. A $d$-bubble (closed or open) in $\cG_{\mf^0}$ is denoted by $\bee^d$, the set of $d$-bubbles by $\cB^d$, and its cardinal $B^d$. The set of vertices and edges of $\bee^d$ are denoted by $\cV_{\bee^d}$ and $\cE_{\bee^d}$ of cardinal  $V_{\bee^d}$ and $E_{\bee^d}$ respectively. We also denote by  $\cF_{\inter; \bee^d}$ and $\cB^p_{\bee^d}$ , the sets of internal faces and $p$-bubbles ($p \leq d$) of $\bee^d$ of cardinal $F_{\inter; \bee^d}$ and $B^p_{\bee^d}$ respectively. 

\begin{definition}[Topological invariant for rank $n$ w-colored 
graph \cite{Avo16}]\label{def:topotens}

Let $\mG(\cV,\cE,\mf^0)$ be a rank n w-colored graph and $\boldsymbol{\alpha}=\{\alpha_k\}_{k=3, \cdots, n}$ some positive rational numbers.
The generalized topological invariant associated with $\mG$
is given by the following  function associated with any of its 
representatives $\cG_{\mf^0}$.
\bea\label{ttopfla} 
&&
\mT_{\mG;\boldsymbol{\alpha}}(x,y,z,s,w,q,t) =\mT_{\cG_{\mf^0};\boldsymbol{\alpha}}(x,y,z,s,w,q,t)= \label{ttopflan}\\
 && \sum_{s \sset  \cG_{\mf^0}}
 (x-1)^{\rk( \cG_{\mf^0})-\rk(s)}y^{n(s)}
z^{\frac{(n-1)(n+2)}{2}k(s)-\gamma_{n;\boldsymbol{\alpha}}(s)}s^{C_\partial(s)} \,  w^{F_{\partial}(s)} q^{E_{\partial}(s)} t^{f(s)}\,, 
\nonumber
\eea
with 
\bea
\gamma_{n;\boldsymbol{\alpha}}(s)=\frac{n(n-1)}{2}(V(s)-E(s)) + (n-1)F_{\inter}(s) - (2+(n-2)\alpha_3)B^3(s) + \\\nonumber\sum_{k=4}^n \Big[(k-1)\alpha_{k-1}-(n-k+1)\alpha_k\Big]B^k(s)
\eea
a negative integer.
\end{definition}

Expanding the definitions of $\rk(s)$ and $n(s)$ in \eqref{ttopfla} yields 
\bea
&&
\mT_{\mG;\boldsymbol{\alpha}}(x,y,z,s,w,q,t) =\mT_{\cG_{\mf^0};\boldsymbol{\alpha}}(x,y,z,s,w,q,t)= \label{ttopflan}\nonumber \\
 && (x-1)^{-k(\cG)}\Big(yz^{\frac{n(n-1)}{2}}\Big)^{-v(\cG)}\sum_{A \sset  \cG_{\mf^0}}
 \Big((x-1)yz^{\frac{(n-1)(n+2)}{2}}\Big)^{k(s)}\Big(yz^{\frac{n(n-1)}{2}}\Big)^{e(s)}z^{(1-n)F_{\inter}(s)}\nonumber\\
&&\times z^{(2+(n-2)\alpha_3)B^3(s)}\Big(\prod_{k=4, \cdots, n}z^{((1-k)\alpha_{k-1}+(n-k+1)\alpha_k)B^k(s)}\Big)s^{C_\partial(s)} \,  w^{F_{\partial}(s)} q^{E_{\partial}(s)} t^{f(s)}.\, 
\eea

Let us introduce the multivariate version of this polynomial 
\begin{definition}[Multivariate form]
The multivariate form associated with \eqref{ttopfla} is defined by:
 \bea
&&
\widetilde\mT_{\mG;\boldsymbol{\alpha}}(x,\{\beta_e\},\{z_i\}_{i=1,2,3},s,w,q,t)= 
\widetilde\mT_{\cG_{\mf^0}}(x,\{\beta_e\},\{z_i\}_{i=1,\cdots,n},z,s,w,q,t) \label{tmulti}\\
&& =   \sum_{A \sset \cG}
 x^{\rk(s)}\Big(\prod_{e\in s}\beta_e\Big)
\Big(\prod_{i=1, \cdots, n}z_i^{B^i}\Big)z^{F_{\inter}(s)}\, 
s^{C_\partial(s)} \,  w^{F_{\partial}(s)} q^{E_{\partial}(s)} t^{f(s)}\,,
\nonumber
\eea
for $\{\beta_e\}_{e\in \cE}$ labelling the edges of the graph $\cG$. 
\end{definition}
\begin{figure}[h]
 \centering
     \begin{minipage}[t]{.8\textwidth}
      \centering
\includegraphics[angle=0, width=7.5cm, height=3.5cm]{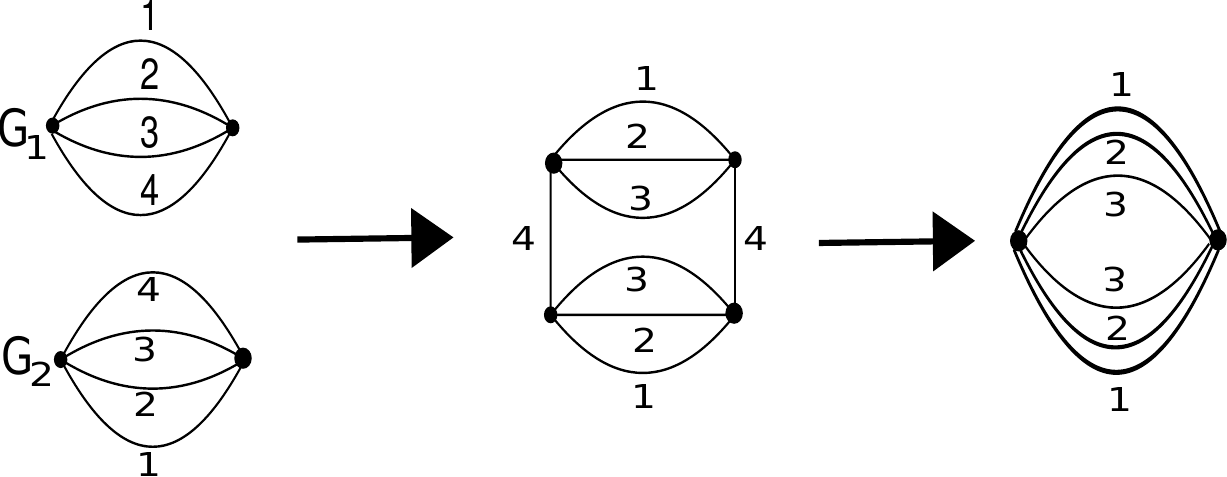}
\caption{ {\small The $2-$sum of $G_1$ and $G_2$ by their edges of color $4$ }}
\label{fig:twosepcolor2}
\end{minipage}
\end{figure}
The operation of $2-$sum introduced earlier can be generalized on the class of graph called weakly colored stranded graph. The only issue is the colors of the edges we want to identify. It is clear that for this operation to be possible, the edges must have the same color. This is illustrated in Figure \ref{fig:twosepcolor2}. 

The construction of  a rank $n$ weakly colored graph $\widehat{\cG}_{\mf^0}$ from the $2-$decomposition $(\cG_{\mf^0},\{H_e\}_{e\in\cE})$  at an edge $e=(u_e,w_e)$ of the template $\cG_{\mf^0}$ will be a direct extension of the one introduced earlier. We can say that the graph $\widehat{\cG}_{\mf^0}$  is obtained  by the identification of the $n-1$ ball (in topological point of view) $m_e$ and $n_e$ of the vertices $u_e$ and $w_e$ in $H_e$ with their correspondence in $\cG_{\mf^0}\vee e$. From this, we can give a partition of $S(H_e)$: the set of states for a given $H_e$. 

Let $S^2(H_e)$ be the set of states of $H_e$ in which $m_e$ and $n_e$ lie in different connected components and different connected components of the boundary graph. We set $S^1(H_e)$ as the set of states of $H_e$ in which $m_e$ and $n_e$ lie in the same connected components and same  internal face or connected component of the boundary graph.
 We work under the condition that $m_e$ and $n_e$ lie in the same connected components if and only if they are in the same internal face or connected component of the boundary graph.
\begin{lemma}\label{lemma:partialfinc}
Let $\hat{s}$ be  a state of the rank $n$ weakly colored stranded graph $\widehat{\cG}_{\mf^0}$. $s_e\in S^1(H_e)\cup S^2(H_e)$, $e\in\cE$ in the decomposition above. Then
\bea\label{eq:faceoxpdern}
k(\hat{s})=\sum_{e\in\cE}k(s_e)+k(s)- |\{s_e\in S^1(H_e)\}|-2|\{s_e\in S^2(H_e)\}|,
\eea
\bea\label{eq:faceopdern}
 &&C_\partial(\hat{s})=\sum_{e\in\cE}C_\partial(s_e)+C_\partial(s)- |\{s_e\in S^1(H_e)\}|-2|\{s_e\in S^2(H_e)\}|,\cr &&F_{\inter}(\hat{s}) = \sum_{e\in\cE}F_{\inter}(s_e)+F_{\inter}(s),
\eea
\bea\label{eq:faceopenc}
 f(\hat{s})=\sum_{e\in\cE}f(s_e)+f(s)- 2|\{s_e\in S^1(H_e)\}|- 4|\{s_e\in S^2(H_e)\}|,
\eea
\bea\label{eq:faceopenss}
 E_\partial(\hat{s})=\sum_{e\in\cE}E_\partial(s_e)+E_\partial(s)- n|\{s_e\in S^1(H_e)\}|-2\times n|\{s_e\in S^2(H_e)\}|,
\eea
\bea\label{eq:faceopen}
 F_\partial(\hat{s})\geq\sum_{e\in\cE}F_\partial(s_e)+F_\partial(s)-n|\{s_e\in S^1(H_e)\}|-2\times n|\{s_e\in S^2(H_e)\}|.
\eea
Furthermore
\bea\label{eq:facecloseds}
B^p(\hat{s}) = \sum_{e\in\cE}B^p(s_e) + B^p(s)-\complement_{n}^{p-1} |\{s_e\in S^1(H_e)\}|-2\complement_{n}^{p-1}|\{s_e\in S^2(H_e)\}|.
\eea
\end{lemma}

\proof
The proof of this lemma can be performed on colored tensor graphs since the w-colored stranded graphs are obtained by a successive contraction of edges in a colored tensor graph. Furthermore this contraction does not modify the boundary or the number of bubbles in a colored tensor graph.

Equations \eqref{eq:faceoxpdern}, \eqref{eq:faceopdern} and \eqref{eq:faceopenc} are direct extension to the relations in Lemma \ref{lemma:partialfincc} .

Let us consider the graph $\partial(s)$.  Each edge of $\partial(s)$ or open face of $s$ corresponds to an edge in $\partial(\hat{s})$.  However, $\partial(\hat{s})$ has extra edges which are the edges of $\partial(s_e)$; $s_e\in S^1(H_e)\cup S^2(H_e)$, $e\in\cE(\cG_1)$. Since $s_e$ is inserted in $s$ by $n$ strands then there are $E(\partial(s_e))-n$ such extra edges. This ends the proof of \eqref{eq:faceopenss}. The proof of \eqref{eq:faceopen} is similar to \eqref{eq:faceopenss}. In this case $s_e$ is inserted in $s$ by at most $3$ faces of the boundary  $\partial(s_e)$.

Let us assume that the graphs $s_e$ are inserted in $s$ by the edges of colors $i$. There is a one to one correspondence between the $p-$bubbles of $s$ and $s_e$ not containing any $i$ and those of $\hat{s}$ not containing the same colors $i$. Furthermore each $s_e\in S^1(H_e)$ share $\complement_{n}^{p-1}$ number of $p-$bubbles with $s$ (these bubbles contain the color of $e$). If  $s_e\in S^2(H_e)$, then the number is $2\complement_{n}^{p-1}$. Hence the number of extra bubbles in $\hat{s}$  are $B^p(s_e)-\complement_{n}^{p-1}$. This ends the proof of \eqref{eq:facecloseds}.

\qed

We consider the following state sums:
\bea\label{eq:eta}
\eta_e^{(1)}(a,b,c,d,f,\{g_p\}) :=\sum_{s\in S^1(H_e)}a^{k(s)-1}b^{e(s)}c^{F_{\inter}(s) }d^{C_\partial(s)-1}f^{E_\partial(s)-n}\Big(\prod_{p=3}^{p=n} g_p^{B^p(s)-\complement_{n}^{p-1}}\Big),
\eea
\bea\label{eq:etaaa}
\eta_e^{(2)}(a,b,c,d,f,\{g_p\}) :=\sum_{s\in S^2(H_e)}a^{k(s)-2}b^{e(s)}c^{F_{\inter}(s) }d^{C_\partial(s)-2}f^{E_\partial(s)-2n}\Big(\prod_{p=3}^{p=n} g_p^{B^p(s)-2\complement_{n}^{p-1}}\Big).
\eea
We can observe in the expressions of $\eta_e^{(1)}$ and $\eta_e^{(2)}$ given in \eqref{eq:eta} and \eqref{eq:etaaa} respectively that their is no variable for the faces of the boundary graph.
The reason of this choice come from the inequality given by \eqref{eq:faceopen}. 

\begin{proposition}\label{prop:hjj}
Let $(\cG_{\mf^0},\{H_e\}_{e\in \cE})$ be a $2-$decomposition of $\widehat{\cG}_{\mf^0}$ and $\eta_e^{(1)}$, $\eta_e^{(2)}$ two functions as introduced above. Then
\bea
Z(\widehat{\cG}_{\mf^0};a,b,c,d,f,\{g_p\}) = &&\sum_{s\in S(\cG)} a^{k(s)}c^{F_{\inter}(s)}d^{C_\partial(s)}f^{E_\partial(s)}\Big(\prod_{p=3}^{p=n} g_p^{B^p(s)}\Big)\times \cr&&\Big(\prod_{e\in s}\eta_e^1\Big)\Big(\prod_{e\notin s}\eta_e^2\Big).
\eea
\end{proposition}

Once again we can make a partition of the set $S(H_e)$ of states $s$ of $H_e$ where $e=(m_e,n_e)$ . We can set $S(H_e)= S^1(H_e)\cup S^2(H_e)$; where $S^1(H_e)$ is the set of states having $m_e$ and $n_e$ in the same connected component and the same component of the boundary graph and $S^2(H_e)$ is the set of states having $m_e$ and $n_e$ in different connected component and different component of the boundary graph. It is clear that the insertion of $e$ in $s\in S^2(H_e)$ decreases the number of connected components of the boundary graph by one, the number of edges of the boundary by $n$ and the number of $p-$bubbles by $\complement_{n}^{p-1}$. Then if $s$ contributes with the term $a^{k(s)}c^{F_{\inter}(s)}d^{C_\partial(s)}f^{E_\partial(s)}\Big(\prod_{p=3}^{p=n} g_p^{B^p(s)}\Big)$, the state $s\cup e$ will contribute with $a^{k(s)-1}c^{F_{\inter}(s)}d^{C_\partial(s)-1}f^{E_\partial(s)-n}\Big(\prod_{p=3}^{p=n} g_p^{B^p(s)-\complement_{n}^{p-1}}\Big)$. 

The insertion of $e$ in a state $s\in S^1(H_e)$ of $H_e$ leads to different possible cases. For example this insertion can increase the number of internal faces from $0$ to $n$. We also have multiple possibilities for the number of connected components of the boundary graph, the number of edges of the boundary graphs and the number of $p-$bubbles. A good analysis of the different possibilities may help to find a more explicit formula than the one given in Proposition \ref{prop:hjj}. For doing this let us make a restriction to $n=3$.

Consider $s\in S^1(H_e)$ and insert $e$ in $s$. The number of $3-$bubbles may decrease from $0$ to $3$; but the number of internal faces may increase from $0$ to $3$. Depending on the number of internal faces we add by the insertion of $e$, let us discuss the different cases for the number of connected components of the boundary graph and the number of edges of the boundary graph.

$\bullet$ Assume that the number of internal faces is stable after the insertion of $e$. The number of connected components of the boundary graph is stable after the insertion of $e$ or may increase by one. The number of edges of the boundary graph will decrease by $3$. Adding the previous discussions about the evolving of the number of $3-$bubbles under the insertion of $e$, if $s$ contribute with the term $$a^{k(s)}c^{F_{\inter}(s)}d^{C_\partial(s)}f^{E_\partial(s)}\Big(\prod_{p=3}^{p=n} g_p^{B^p(s)}\Big),$$ in the polynomial, then the state $s\cup e$ will contribute with $a^{k(s)}c^{F_{\inter}(s)}d^{C_\partial(s)+\beta(s)}f^{E_\partial(s)-3}g_3^{B^3(s)-\gamma(s)}$ for $\beta(s)=0,1$ and $\gamma(s)=0,1,2,3$.

$\bullet$ Assume that the number of internal faces increases by one or two after the insertion of $e$. If the number of internal faces increases by one or two after the insertion of $e$, then number of edges of the boundary graph will decrease by two or one respectively. 
It is also clear that the number of connected components of the boundary graph is stable or may increase by one after the insertion of $e$. If $s$ contribute with the term $a^{k(s)}c^{F_{\inter}(s)}d^{C_\partial(s)}f^{E_\partial(s)}\Big(\prod_{p=3}^{p=n} g_p^{B^p(s)}\Big)$, in the polynomial, then the state $s\cup e$ will contribute with $$a^{k(s)}c^{F_{\inter}(s)+\alpha(s)}d^{C_\partial(s)+\beta(s)}f^{E_\partial(s)-3+\alpha(s)}g_3^{B^3(s)-\gamma(s)},$$ for $\alpha(s)=1,2$, $\beta(s)=0,1$ and $\gamma(s)=0,1,2,3$.

$\bullet$ Assume that the number of internal faces increases by three after the insertion of $e$. The number of connected components of the boundary graph is stable after the insertion of $e$ or decreases by one. The number of edges of the boundary graphs is stable after the insertion. If $s$ contributes with the term $a^{k(s)}c^{F_{\inter}(s)}d^{C_\partial(s)}f^{E_\partial(s)}\Big(\prod_{p=3}^{p=n} g_p^{B^p(s)}\Big)$, in the polynomial, then the state $s\cup e$ will contribute with $a^{k(s)}c^{F_{\inter}(s)+3}d^{C_\partial(s)+\beta(s)}f^{E_\partial(s)}g_3^{B^3(s)-\gamma(s)}$ for $\beta(s)=-1,0$ and $\gamma(s)=0,1,2,3$. 

We can summarized all the different cases in this way: if $s$ contributes with the term $$a^{k(s)}c^{F_{\inter}(s)}d^{C_\partial(s)}f^{E_\partial(s)}\Big(\prod_{p=3}^{p=n} g_p^{B^p(s)}\Big),$$ in the polynomial, then the state $s\cup e$ will contribute with $$a^{k(s)}c^{F_{\inter}(s)+\alpha(s)}d^{C_\partial(s)+\beta(s)}f^{E_\partial(s)-3+\alpha(s)}g_3^{B^3(s)-\gamma(s)},$$ for $\alpha(s)=0,1,2,3$, $\beta(s)=-1,0,1$ and $\gamma(s)=0,1,2,3$. 

Consider the following sum 
\bea
Z_{A_e}(a,b,c,d,f,g_3)&=&\sum_{s\in S(A_e)}a^{k(s)}c^{F_{\inter}(s)}d^{C_\partial(s)}f^{E_\partial(s)}g_3^{B^3(s)},
\eea
which is also equal to
\bea\label{uhg}
Z_{A_e}(a,b,c,d,f,g_3)&=&\sum_{s\in S^{1}(H_e)}a^{k(s)}c^{F_{\inter}(s)}d^{C_\partial(s)}f^{E_\partial(s)}g_3^{B^3(s)}\cr &+&x_e\sum\limits_{\substack{s\in S^{1}(H_e)\\ \beta(s)=-1,0,1 \\ \gamma(s),\alpha(s)=0,1,2,3}}a^{k(s)}c^{F_{\inter}(s)+\alpha(s)}d^{C_\partial(s)+\beta(s)}f^{E_\partial(s)-3+\alpha(s)} g_3^{B^3(s)-\gamma(s)}\cr &+&(1+a^{-1}x_ed^{-1}f^{-3}g_3^{-3})\sum_{s\in S^{2}(H_e)}a^{k(s)}c^{F_{\inter}(s)}d^{C_\partial(s)}f^{E_\partial(s)} g_3^{B^3(s)}.
\eea
Furthermore 
\bea\label{dfgh}
Z_{A_e}(a,b,c,d,l)=Z_{H_e}(a,b,c,d,l)+x_eZ_{A_e/e}(a,b,c,d,l).
\eea
A reformulation of the equation \eqref{uhg} in a great number of terms together with an identification with \eqref{dfgh} will help us to find a theorem similar to Theorem \ref{theo:prod}. The case where the points $m_e$ and $n_e$ may  belong to the same connected component but different boundary components is not studied in this work and remains an open question for future investigations. 
\begin{center}
{\bf Acknowledgements}
\end{center}

{\footnotesize The author wishes to express his gratitude to Joseph Ben Geloun  for his assistance in the presentation of this paper. The author thanks the Max-Planck Institute, Albert Einstein Institute for its hospitality.}

\vspi

\vspace{0.5cm}

\end{document}